\documentclass[10.5pt]{article}
\usepackage{amsmath}
\usepackage{amsfonts}
\usepackage{color}

\def\sqr#1#2{{\vcenter{\vbox{\hrule height.#2pt
              \hbox{\vrule width.#2pt height#1pt \kern#1pt \vrule width.#2pt}
              \hrule height.#2pt}}}}
\def\signed #1{{\unskip\nobreak\hfil\penalty50
              \hskip2em\hbox{}\nobreak\hfil#1
              \parfillskip=0pt \finalhyphendemerits=0 \par}}
\def\endpf{\signed {$\sqr69$}}

\def\dbR{{\mathop{\rm l\negthinspace R}}}
\def\dbC{{\mathop{\rm l\negthinspace\negthinspace\negthinspace C}}}

\def\5n{\negthinspace \negthinspace \negthinspace \negthinspace \negthinspace }
\def\4n{\negthinspace \negthinspace \negthinspace \negthinspace }
\def\3n{\negthinspace \negthinspace \negthinspace }
\def\2n{\negthinspace \negthinspace }
\def\1n{\negthinspace }

\def\dbC{\mathbb{C}}

\def\dbE{\mathbb{E}}
\def\dbF{\mathbb{F}}

\def\dbP{\mathbb{P}}

\def\dbR{\mathbb{R}}

\def\={\buildrel \triangle \over =}

\def\ds{\displaystyle}

\def\ns{\noalign{\ss}}

\def\d{\delta}

\def\l{\lambda}
\def\m{\mu}

\def\si{\sigma}
\def\t{\tau}

\def\o{\omega}

\def\R{{\bf R}}
\def\G{\Gamma}
\def\D{\Delta}
\def\Th{\Theta}
\def\L{\Lambda}

\def\F{\Phi}
\def\O{\Omega}

\def\cF{{\cal F}}

\def\cN{{\cal N}}

\def\cP{{\cal P}}

\def\cR{{\cal R}}
\def\cS{{\cal S}}

\def\cU{{\cal U}}

\def\cX{{\cal X}}

\def\A{\bar{A}}
\def\B{\bar{B}}
\def\C{\bar{C}}
\def\D{\bar{D}}
\def\F{\bar{F}}
\def\K{\bar{K}}
\def\L{\bar{L}}
\def\Q{\bar{Q}}
\def\R{\bar{R}}
\def\S{\bar{S}}
\def\T{\bar{T}}
\def\U{\bar{U}}
\def\V{\bar{V}}
\def\W{\bar{W}}
\def\X{\bar{X}}
\def\Y{\bar{Y}}

\def\no{\noindent}

\def\ss{\smallskip}
\def\ms{\medskip}

\def\q{\quad}
\def\qq{\qquad}
\def\hb{\hbox}

\def\limsup{\mathop{\overline{\rm lim}}}

\def\Ra{\mathop{\Rightarrow}}

\def\lan{\mathop{\langle}}
\def\ran{\mathop{\rangle}}

\def\wt{\widetilde}
\def\cd{\cdot}
\def\cds{\cdots}

\def\deq{\mathop{\buildrel\Delta\over=}}

\def\var{{\mathop{\rm var}\,}}
\def\({\Big (}
\def\){\Big )}
\def\[{\Big[}
\def\]{\Big]}
\def\bde{\begin{definition}}
\def\ede{\end{definition}}
\def\be{\begin{equation}}
\def\bel{\begin{equation}\label}
\def\ee{\end{equation}}
\def\bex{\begin{example}}
\def\eex{\end{example}}
\def\bt{\begin{theorem}}
\def\et{\end{theorem}}
\def\bc{\begin{corollary}}
\def\ec{\end{corollary}}
\def\bl{\begin{lemma}}
\def\el{\end{lemma}}
\def\bp{\begin{proposition}}
\def\ep{\end{proposition}}
\def\bas{\begin{assumption}}
\def\eas{\end{assumption}}
\def\br{\begin{remark}}
\def\er{\end{remark}}
\def\ba{\begin{array}}
\def\ea{\end{array}}
\def\ed{\end{document}}

\def\square#1{\vbox{\hrule\hbox{\vrule height#1%
     \kern#1\vrule}\hrule}}
\def\rectangle#1#2{\vbox{\hrule\hbox{\vrule height#1%
     \kern#2\vrule}\hrule}}

\font\tenbb=msbm10
\font\sevenbb=msbm7
\font\fivebb=msbm5

\newfam\bbfam
\scriptscriptfont\bbfam=\fivebb
\textfont\bbfam=\tenbb
\scriptfont\bbfam=\sevenbb

\newtheorem{theorem}{Theorem}[section]
\newtheorem{lemma}[theorem]{Lemma}
\newtheorem{remark}[theorem]{Remark}
\newtheorem{example}[theorem]{Example}
\newtheorem{corollary}[theorem]{Corollary}
\newtheorem{definition}[theorem]{Definition}
\newtheorem{proposition}[theorem]{Proposition}
\newtheorem{assumption}[theorem]{Assumption}

\makeatletter
   
   \@addtoreset{equation}{section}
\makeatother

\begin{document}

\title{\bf A Linear-Quadratic Optimal Control Problem for Mean-Field Stochastic
Differential Equations in Infinite Horizon\footnote{This work is
supported in part by RGC Grants GRF521610 and GRF501010, and NSF
Grant DMS-1007514.}}

\author{Jianhui Huang$^a$, Xun Li$^a$,~~and~~Jiongmin Yong$^b$ \\
{\footnotesize\textsl{$^a$Department of Applied Mathematics, The Hong Kong Polytechnic University, Hong Kong, China}} \\
{\footnotesize\textsl{$^b$Department of Mathematics, University of Central Florida, Orlando, FL 32816, USA}}}

\maketitle

\begin{abstract}
A linear-quadratic (LQ, for short) optimal control problem is considered
for mean-field stochastic differential equations with constant
coefficients in an infinite horizon. The stabilizability of the
control system is studied followed by the discussion of the
well-posedness of the LQ problem. The optimal control can be
expressed as a linear state feedback involving the state and its
mean, through the solutions of two algebraic Riccati equations. The
solvability of such kind of Riccati equations is investigated by
means of semi-definite programming method.
\end{abstract}

\ms


\ms


\section{Introduction.}

Let $(\O,\cF,\dbP,\dbF)$ be a complete filtered probability space,
on which a one-dimensional standard Brownian motion $W(\cd)$ is
defined with $\dbF\equiv\{\cF_t\}_{t\ge0}$ being its natural
filtration augmented by all the $\dbP$-null sets. Consider the
following controlled linear stochastic differential equation (SDE,
for short) in $\dbR^n$:
\bel{MF-state1}\left\{\3n\ba{ll}
\ns\ds dX(t)=\Big\{AX(t)+\bar A\dbE[X(t)]+Bu(t)+\bar B\dbE[u(t)]\Big\}dt \\
\ns\ds\qq\qq+\Big\{CX(t)+\bar C\dbE[X(t)]+Du(t)+\bar D\dbE[u(t)]\Big\}dW(t),\q t\ge0, \\
\ns\ds X(0)=x,\ea\right.\ee
where $A,\bar A,C,\bar C\in\dbR^{n\times n}$ and $B,\bar B,D,\bar D\in\dbR^{n\times m}$
are given (deterministic) matrices. In the
above, $X(\cd)$, valued in $\dbR^n$, is called the {\it state process},
and $u(\cd)$, valued in $\dbR^m$, is called a {\it control process}.

\ms

Different from classical controlled linear SDEs, the terms
$\dbE[X(\cd)]$ and $\dbE[u(\cd)]$ appear in the equation. We call
(\ref{MF-state1}) a controlled mean-field (forward) SDE (MF-FSDE,
for short). Historically, a special case of MF-FSDE, called
McKean--Vlasov SDE, was suggested by Kac \cite{Kac 1956} in 1956 as
a stochastic toy model for the Vlasov type kinetic equation of
plasma and the rigorous study of which was initiated by McKean
\cite{McKean 1966} in 1966. Since then, such kind of equations were
studied by many authors, see, for examples, Dawson \cite{Dawson
1983}, Dawson--G\"artner \cite{Dawson-Gartner 1987}, G\'artner
\cite{Gartner 1988}, Scheutzow \cite{Scheutzow 1987}, Graham
\cite{Graham 1992}, Chan \cite{Chan 1994}, Chiang \cite{Chiang
1994}, Ahmed--Ding \cite{Ahmed-Ding 1995}, and the references cited
therein. For the relevant works of recent years, see, for examples,
Veretennikov \cite{Veretennikov 2003}, Huang--Malham\'e--Caines
\cite{Huang-Malhame-Caines 2006}, Buckdahn-Djehiche-Li-Peng
\cite{Buckdahn-Djehiche-Li-Peng 2009}, Buckdahn-Li-Peng
\cite{Buckdahn-Li-Peng 2009}, Borkar--Kumar \cite{Borkar-Kumar
2010}, Crisan--Xiong \cite{Crisan-Xiong 2010}, Kotelenez--Kurtz
\cite{Kotelenez-Kurtz 2010}, and so on. Control problems of
McKean--Vlasov equation or MF-FSDEs were investigated by Ahmed--Ding
\cite{Ahmed-Ding 2001}, Ahmed \cite{Ahmed 2007},
Buckdahn-Djehiche-Li \cite{Buckdahn-Djehiche-Li 2010}),
Park--Balasubramaniam--Kang \cite{Park-Balasubramaniam-Kang 2008},
Andersson--Djehiche \cite{Andersson-Djehche 2011},
Meyer-Brandis--Oksendal--Zhou \cite{Meyer-Brandis-Oksandal-Zhou
2011}, and so on. In Yong \cite{Yong 2011}, a linear-quadratic (LQ,
for short) problem was introduced and investigated for MF-FSDEs in
finite horizons. Some interesting motivation was given in \cite{Yong
2011} for the control problem with $\dbE[X(\cd)]$ and $\dbE[u(\cd)]$
being included in the cost functional. This paper can be regarded as
a continuation of \cite{Yong 2011}, for LQ problem of MF-FSDEs in an
infinite horizon.

\ms

We introduce the following:
%
$$\left\{\2n\ba{ll}
\ns\ds\cU[0,T]=\Big\{u:[0,\infty)\times\O\to\dbR^m\bigm|u(\cd)\hb{ is $\dbF$-adapted, } \dbE\int_0^T|u(s)|^2ds<\infty\Big\},\qq\forall T>0, \\
\ns\ds\cU_{loc}[0,\infty)=\bigcup_{T>0}\cU[0,T], \\
\ns\ds\cU[0,\infty)=\Big\{u(\cd)\in\cU_{loc}[0,\infty)\bigm|\dbE\int_0^\infty|u(s)|^2ds<\infty\Big\}.\ea\right.$$
Any $u(\cd)\in\cU_{loc}[0,\infty)$ is called a {\it control process}
and any $u(\cd)\in\cU[0,\infty)$ is called a {\it feasible control process}. Likewise, we define
%
$$\left\{\2n\ba{ll}
\ns\ds\cX[0,T]=\left\{X:[0,\infty)\times\O\to\dbR^n\left\vert\ba{l} X(\cd)\hb{ is $\dbF$-adapted, }t\mapsto X(t,\o)\hb{ is continuous}, \\
\ns\ds\dbE\big[\sup_{t\in[0,T]}|X(t)|^2\big]<\infty\ea\right.\right\},\qq T>0, \\
\ns\ds\cX_{loc}[0,\infty)=\bigcup_{T>0}\cX[0,T], \\
\ns\ds\cX[0,\infty)=\Big\{X(\cd)\in\cX_{loc}[0,\infty)\bigm|\dbE\int_0^\infty|X(t)|^2dt<\infty\Big\}.\ea\right.$$
Any element in $\cX_{loc}[0,\infty)$ is called a {\it state process}. It is not hard to see that
%
$$\left\{\2n\ba{ll}
\ns\ds\cU[0,\infty)\subseteq\cU_{loc}[0,\infty), & \q \cU[0,\infty)\ne\cU_{loc}[0,\infty), \\
\ns\ds\cX[0,\infty)\subseteq\cX_{loc}[0,\infty), & \q \cX[0,\infty)\ne\cX_{loc}[0,\infty).\ea\right.$$
By a standard argument using contraction mapping theorem, one can
show that for any $(x,u(\cd))\in\dbR^n\times\cU_{loc}[0,\infty)$,
(\ref{MF-state1}) admits a unique solution
$X(\cd)=X(\cd\,;x,u(\cd))\in\cX_{loc}[0,\infty)$. Next, we let
$Q,\bar Q\in\cS^n$ and $R,\bar R\in\cS^m$, where $\cS^k$ is the set
of all symmetric matrices of order $(k\times k)$, and introduce the
following cost functional:
\bel{J(x,u)}\ba{ll}
\ns\ds J(x;u(\cd))=\dbE\int_0^\infty\Big\{\lan QX(s),X(s)\ran+\lan\bar Q\dbE[X(s)],\dbE[X(s)]\ran \\
\ns\ds\qq\qq\qq\qq\qq+\lan Ru(s),u(s)\ran+\lan\bar
R\dbE[u(s)],\dbE[u(s)]\ran\Big\}ds,\ea\ee
where $X(\cd)=X(\cd\,;x,u(\cd))$ on the right hand side of the
above. Note that in general, for
$(x,u(\cd))\in\dbR^n\times\cU[0,\infty)$, the solution $X(\cd)\equiv
X(\cd\,;x,u(\cd))$ of (\ref{MF-state1}) might just be in
$\cX_{loc}[0,\infty)$ and the above cost functional $J(x;u(\cd))$
might not be defined. Therefore, we introduce the following:
$$\cU_{ad}[0,\infty)=\Big\{u(\cd)\in\cU[0,\infty)\bigm|J(x;u(\cd))\hb{ is defined, }\forall x\in\dbR^n\Big\}.$$
Any element $u(\cd)\in\cU_{ad}[0,\infty)$ is called an {\it
admissible control process} and the corresponding $X(\cd)\equiv
X(\cd\,;x,u(\cd))$ is called an {\it admissible state process}. We
see that the structure of $\cU_{ad}[0,\infty)$ is very complicated,
since it involves not only the state equation, but also the cost
functional. Some better description of $\cU_{ad}[0,\infty)$ will be
given a little later, under proper conditions. Our optimal control
problem can be stated as follows:

\ms

\bf Problem (MF-LQ). \rm For given $x\in\dbR^n$, find a
$u_*(\cd)\in\cU_{ad}[0,\infty)$ such that
%
$$J(x;u_*(\cd))=\inf_{u(\cd)\in\cU_{ad}[0,\infty)}J(x;u(\cd))\equiv V(x).$$

\ms

Any $u_*(\cd)\in\cU_{ad}[0,\infty)$ satisfying the above is called
an {\it optimal control process} and the corresponding state process
$X_*(\cd)\equiv X(\cd\,;x,u_*(\cd))$ is called an {\it optimal
state process}; the pair $(X_*(\cd),u_*(\cd))$ is called an {\it
optimal pair}. The function $V(\cd)$ is called the {\it value
function} of Problem (MF-LQ).

\ms

It is not hard to see that in order Problem (MF-LQ) to make sense,
we need $\cU_{ad}[0,\infty)$ to be nonempty, at least. To achieve
this, we will carefully discuss various stabilizability (for which
both the state equation and the cost functional are involved) of the
controlled MF-FSDE (\ref{MF-state1}), which are interestingly
different from the classic ones, due the the appearance of the terms
$\dbE[X(\cd)]$ and $\dbE[u(\cd)]$. Once the set $\cU_{ad}[0,\infty)$
of admissible controls is nonempty, under some standard assumptions,
we are able to show that the optimal control uniquely exists. Then
inspired by the results of \cite{Yong 2011}, we obtain a system of
algebraic Riccati equations (AREs, for short), whose solutions will
lead us to the state feedback representation of the optimal control.
The existence of the solutions to the derived ARE system is
established under some reasonable conditions. Our results recovers
relevant ones for the classic linear-quadratic optimal controls of
SDEs.

\ms

The rest of the paper is organized as follows. Section 2 collects
some preliminary results concerning the state equation. In Sections
3 and 4, the stability and the stabilizability of the state equation
are discussed. In Section 5, Problem (MF-LQ) is solved by means of
AREs. In Section 6, the solvability of AREs is discussed by linear
matrix inequalities (LMIs, for short). A couple of numerical
examples are presented in Section 7. Finally, some supporting
results for Section 6 are listed in the Appendix.

\section{Preliminary Results}

In this section, we present some preliminary results. First of all,
let us consider the following result, whose proof follows a standard
argument using contraction mapping theorem, together with It\^o's
formula.

\ms

\bp\label{prop2.1} \sl
For any $x\in\dbR^n$ and
$u(\cd)\in\cU[0,\infty)$, there exists a unique $X(\cd)\equiv
X(\cd\,;x,u(\cd))$ solving $(\ref{MF-state1})$. Moreover,
%
$$\dbE\big[\sup_{t\in[0,T]}|X(t)|^2\big]\le L_T\Big\{|x|^2+\dbE\int_0^T|u(t)|^2dt\Big\},\qq\forall T>0.$$
where $L_T>0$ is a constant depending on $T$, and independent of
$(x,u(\cd))$.
\ep

\ms

\rm

For later purposes, we make some calculations. Let
$X(\cd)=X(\cd\,;x,u(\cd))$ be the solution of (\ref{MF-state1}). For
any deterministic differentiable function $P(\cd)$ valued in
$\cS^n$, by It\^o's formula, we have
$$\ba{ll}
\ns\ds d\lan P(t)X(t),X(t)\ran \\
\ns\ds=\Big\{\1n\lan\dot P(t)X(t),X(t)\ran+2\lan P(t)\big\{AX(t)+\bar A\dbE[X(t)]+Bu(t)+\bar B\dbE[u(t)]\big\},X(t)\ran \\
\ns\ds\q+\1n\lan\1n P(t)\big\{CX(t)\1n+\1n\bar C\dbE[X(t)]\1n+\1n Du(t)\1n+\1n\bar D\dbE[u(t)]\big\},CX(t)\1n+\1n\bar C\dbE[X(t)]\1n+\1n Du(t)\1n+\1n\bar D\dbE[u(t)]\ran\2n\Big\}dt\2n+\1n\{\cds\}dW(t)\\
\ns\ds=\Big\{\lan\dot P(t)X(t),X(t)\ran+2\lan P(t)\big\{AX(t)\1n+\1n Bu(t)\big\},X(t)\ran\1n+\1n\lan P(t)\big\{CX(t)\1n+\1n Du(t)\big\},CX(t)\1n+\1n Du(t)\ran\\
\ns\ds\q+2\lan P(t)\big\{\bar A\dbE[X(t)]+\bar B\dbE[u(t)]\big\}, X(t)\ran\1n+2\lan P(t)\big\{CX(t)+Du(t)\big\},\bar C\dbE[X(t)]+\bar D\dbE[u(t)]\ran\\
\ns\ds\q+\lan P(t)\big\{\bar C\dbE[X(t)]+\bar D\dbE[u(t)]\big\},
\bar C\dbE[X(t)]+\bar D\dbE[u(t)]\ran\1n\Big\}dt+\{\cds\}dW(t).\ea$$
Thus,
\bel{PX,X}\ba{ll}
\ns\ds\dbE\lan P(t)X(t),X(t)\ran=\lan P(0)x,x\ran+\dbE\int_0^t\Big\{\lan\big[\dot P(s)+P(s)A+A^TP(s)+C^TP(s)C\big]X(s),X(s)\ran \\
\ns\ds\qq\qq\qq\qq+2\lan u(s),(B^T\1n P(s)+D^T\1n P(s)C\big]X(s)\ran+\lan D^T\1n P(s)Du(s),u(s)\ran \\
\ns\ds\qq\qq\qq\qq+\lan\big[P(s)\bar A+\bar A^T\1n P(s)+\bar C^T\1n P(s)\bar C+\bar C^T\1n P(s)C+C^T\1n P(s)\bar C\big]\dbE[X(s)],\dbE[X(s)]\ran \\
\ns\ds\qq\qq\qq\qq+2\lan\dbE[u(s)],\big[\bar B^TP(s)+\bar D^TP(s)\bar C+\bar D^TP(s)C+D^TP(s)\bar C\big]\dbE[X(s)]\ran \\
\ns\ds\qq\qq\qq\qq+\lan\big[\bar D^TP(s)\bar D+\bar
D^TP(s)D+D^TP(s)\bar D\big]\dbE[u(s)],\dbE[u(s)]\ran\Big\}ds.\ea\ee
Also,
\bel{PEX,EX}\ba{ll}
\ns\ds\lan P(t)\dbE[X(t)],\dbE[X(t)]\ran=\lan P(0)x,x\ran\1n+\2n\int_0^t\2n\Big\{\1n\lan\1n\big[\dot P(s)\1n+\1n P(s)(A\1n+\1n\bar A)\1n+\1n(A\1n+\1n\bar A)^T\1n P(s)\big]\dbE[X(s)],\dbE[X(s)]\ran \\
\ns\ds\qq\qq\qq\qq\qq\qq\qq+2\lan\dbE[u(s)],(B\1n+\1n\bar B)^T\1n P(s)\dbE[X(s)]\ran\2n\Big\}ds.\ea\ee
Combining (\ref{PX,X}) and (\ref{PEX,EX}), we obtain
%
$$\3n\3n\ba{ll}
\ns\ds\dbE\lan P(t)\big\{X(t)-\dbE[X(t)]\big\},X(t)-\dbE[X(t)]\ran=\dbE\lan P(t)X(t),X(t)\ran-\lan P(t)\dbE[X(t)],\dbE[X(t)]\ran \\
\ns\ds=\dbE\2n\int_0^t\2n\Big\{\2n\lan\big[\dot P(s)\1n+\1n P(s)A\1n+A^T\1n P(s)\1n+\1n C^T\1n P(s)C\big]X(s),X(s)\ran+2\lan u(s),\big[B^T\1n P(s)\1n+\1n D^T\1n P(s)C\big]X(s)\ran \\
\ns\ds\qq +\lan D^T\1n P(s)Du(s),u(s)\ran\1n+\lan\1n\big[P(s)\bar A\1n+\1n\bar A^T\1n P(s)\1n+\1n\bar C^T\1n P(s)\bar C\2n+\1n\bar C^T\1n P(s)C\1n+\1n C^T\1n P(s)\bar C\big]\dbE[X(s)],\dbE[X(s)]\ran\\
\ns\ds\qq +2\lan\dbE[u(s)],\big[\bar B^T\1n P(s)\1n+\1n\bar D^T\1n P(s)\bar C\1n+\1n\bar D^T\1n P(s)C\1n+\1n D^TP(s)\bar C\big]\dbE[X(s)]\ran \\
\ns\ds\qq +\lan\big[\bar D^T\1n P(s)\bar D\1n+\1n\bar D^T\1n P(s)D\1n+\1n D^T\1n P(s)\bar D\big]\dbE[u(s)],\dbE[u(s)]\ran \\
\ns\ds\qq -\lan\big[\dot P(s)+P(s)(A+\bar A)+(A+\bar A)^TP(s)\big]\dbE[X(s)],\dbE[X(s)]\ran-2\lan\dbE[u(s)],(B+\bar B)^TP(s)\dbE[X(s)]\ran\Big\}ds \\
\ns\ds\1n=\1n\dbE\2n\int_0^t\2n\Big\{\lan\big[\dot P(s)\1n+\1n P(s)A\1n+\1n A^T\1n P(s)\1n+\1n C^T\1n P(s)C\big]\big\{X(s)-\dbE[X(s)]\big\},X(s)-\dbE[X(s)]\ran \\
\ns\ds\qq +2\lan u(s)\1n-\1n\dbE[u(s)],\big[B^T\2n P(s)\2n+\2n D^T\2n P(s)C\big]\big\{X(s)\1n-\1n\dbE[X(s)]\big\}\1n\ran \\
\ns\ds\qq +\1n\lan\1n D^T\2n P(s)D\big\{u(s)\1n-\1n\dbE[u(s)]\big\},u(s)\1n-\1n\dbE[X(s)]\ran \\
\ns\ds\qq +\1n\lan P(s)\big\{(C\2n+\1n\bar C)\dbE[X(s)]+(D\1n+\1n\bar D)\dbE[u(s)]\big\},(C\2n+\1n\bar C)\dbE[X(s)]+(D\1n+\1n\bar D)\dbE[u(s)]\ran\Big\}ds.\ea$$
In the case that $P(t)\equiv P\in\cS^n$, we have
\bel{PX,X0}\ba{ll}
\ns\ds\dbE\lan PX(t),X(t)\ran = \ds\lan Px,x\ran+\dbE\int_0^t\Big\{\lan(PA+A^TP+C^TPC)X(s),X(s)\ran \\
\ns\ds\qq\qq\qq +2\lan u(s),(B^T\1n P\1n+\1n D^T\1n PC)X(s)\ran\1n+\1n\lan D^T\1n PDu(s),u(s)\ran \\
\ns\ds\qq\qq\qq +\lan(P\bar A\1n+\1n\bar A^T\1n P\1n+\1n\bar C^T\1n P\bar C\1n+\1n\bar C^T\1n PC\1n+\1n C^T\1n P\bar C)\dbE[X(s)],\dbE[X(s)]\ran \\
\ns\ds\qq\qq\qq +2\lan\dbE[u(s)],(\bar B^TP+\bar D^TP\bar C+\bar D^TPC+D^TP\bar C)\dbE[X(s)]\ran \\
\ns\ds\qq\qq\qq +\lan(\bar D^TP\bar D+\bar D^TPD+D^TP\bar D)\dbE[u(s)],\dbE[u(s)]\ran\Big\}ds,\ea\ee
%
%
$$\ba{ll}
\ns\ds\lan P\dbE[X(t)],\dbE[X(t)]\ran = \ds\lan Px,x\ran+\2n\int_0^t\2n\Big\{\lan\big[P(A\1n+\1n\bar A)\1n+\1n(A\1n+\1n\bar A)^T\1n P\big]\dbE[X(s)],\dbE[X(s)]\ran \\
\ns\ds\qq\qq\qq\qq\qq +2\lan\dbE[u(s)],(B\1n+\1n\bar B)^T\1n P\dbE[X(s)]\ran\2n\Big\}ds\ea$$
and
\bel{P(X-EX),X-EX0}\ba{ll}
\ns\ds\dbE\lan P\big\{X(t)-\dbE[X(t)]\big\},X(t)-\dbE[X(t)]\ran \\
\ns\ds=\dbE\int_0^t\Big\{\lan(PA+A^TP+C^TPC)\big\{X(s)-\dbE[X(s)]\big\},X(s)-\dbE[X(s)]\ran \\
\ns\ds\q+2\1n\lan u(s)\1n-\1n\dbE[u(s)],(B^T\1n P\2n+\1n D^T\1n PC)\big\{X(s)\1n-\1n\dbE[X(s)]\big\}\ran\1n+\1n\lan D^T\1n PD\big\{u(s)\1n-\1n\dbE[u(s)]\big\},u(s)\1n-\1n\dbE[X(s)]\ran \\
\ns\ds\q+\1n\lan P\big\{(C\2n+\1n\bar C)\dbE[X(s)]+(D\1n+\1n\bar D)\dbE[u(s)]\big\},(C\2n+\1n\bar C)\dbE[X(s)]+(D\1n+\1n\bar D)\dbE[u(s)]\ran\Big\}ds.\ea\ee
The above will be useful later.

\ms

Now, let us look at the cost functional. We observe that the cost
functional $J(x;u(\cd))$ defined by (\ref{J(x,u)}) can also be written as
%
$$\ba{ll}
\ns\ds J(x;u(\cd)) = \ds\dbE\int_0^\infty\Big\{\lan Q\big\{X(t)-\dbE[X(t)]\big\},X(t)-\dbE[X(t)]\ran+\lan(Q+\bar Q)\dbE[X(t)],\dbE[X(t)]\ran \\
\ns\ds\qq\qq\qq +\lan R\big\{u(t)-\dbE[u(t)]\big\},u(t)-\dbE[u(t)]\ran+\lan(R+\bar R)\dbE[u(t)],\dbE[u(t)]\ran\Big\}dt.\ea$$
In what follows, when the dimension of a matrix, say, $Q$ is clear
from the context, we write $Q\ge0$ for $Q\in\cS^n$ being positive
semi-definite and write $Q>0$ for $Q\in\cS^n$ being positive
definite. We now introduce the following assumption concerning the
weighting matrices $Q,\bar Q,R,\bar R$ in the cost functional.

\ms

{\bf(J)} The matrices $Q,\bar Q\in\cS^n$ and $R,\bar R\in\cS^m$ satisfy the following:
%
$$Q,\,Q+\bar Q\ge0,\qq R,\,R+\bar R>0.$$

\ms

Note that in (J), we do not have direct assumption on $\bar Q$ and
$\bar R$, they do not have to be positive (semi-) definite, and
actually, they could even be negative definite. Under (J), we see
that $u(\cd)\in\cU_{ad}[0,\infty)$ if and only if for any
$x\in\dbR^n$, the corresponding state process $X(\cd)\equiv
X(\cd\,;x,u(\cd))$ satisfies
\bel{QQ-stable}
\dbE\int_0^\infty\(|Q^{1\over2}\big\{X(t)-\dbE[X(t)]\big\}|^2+|(Q+\bar Q)^{1\over2}\dbE[X(t)]|^2\)dt<\infty. \ee
Since $Q$ and/or $(Q+\bar Q)$ might be degenerate, when
$u(\cd)\in\cU_{ad}[0,\infty)$, we might not have $X(\cd)\equiv
X(\cd,;x,u(\cd))\in\cX[0,\infty)$. The following is a little
stronger assumption than (J).

\ms

{\bf(J)$'$} The matrices $Q,\bar Q\in\cS^n$ and $R,\bar R\in\cS^m$
satisfy the following:
%
$$Q,\,Q+\bar Q>0,\qq R,\,R+\bar R>0.$$

\ms

Clearly, if (J)$'$ holds, then $u(\cd)\in\cU_{ad}[0,\infty)$ if and
only if for any $x\in\dbR^n$, $X(\cd\,;x,u(\cd))\in\cX[0,\infty)$.

\section{Stability}

Now, let us return to state equation (\ref{MF-state1}). We know that
cost functional $J(x;u(\cd))$ is well-defined on
$\dbR^n\times\cU_{ad}[0,\infty)$, and unlike $\cU[0,\infty)$, the
structure of $\cU_{ad}[0,\infty)$ seems to be complicated since it
involves the state equation and the cost functional. Further, the
following example shows that $\cU_{ad}[0,\infty)$ could even be
empty, which leads to that Problem (MF-LQ) is meaningless.

\ms

\bex\label{ex3.1}
\rm Consider one-dimensional controlled system:
$$dX(t)=X(t)dt+\big\{\dbE[X(t)]+u(t)\big\}dW(t),\qq t\ge0,$$
with cost functional
$$J(x;u(\cd))=\dbE\int_0^\infty|X(t)|^2dt.$$
Clearly,
$$d\dbE[X(t)]=\dbE[X(t)]dt,\qq t\ge0,$$
which implies
$$\dbE[X(t)]=xe^t,\qq t\ge0.$$
Then
$$dX(t)=X(t)dt+\big[xe^t+u(t)\big]dW(t),\qq t\ge0.$$
Hence,
%
%
$$X(t)=xe^t+\int_0^te^{t-s}\big[xe^s+u(s)\big]dW(s)=e^t\Big\{x+\int_0^t\[x+e^{-s}u(s)\]dW(s)\Big\}, \qq t\ge0,$$
and as long as $x\ne0$ or $u(\cd)\ne0$,
%
%
$$J(x;u(\cd))=\dbE\int_0^\infty|X(t)|^2dt=\int_0^\infty e^{2t}\Big\{x^2+\int_0^t\big[x+e^{-s}u(s)\big]^2ds\Big\}dt=\infty.$$
Therefore, in this case, $\cU_{ad}[0,\infty)=\emptyset$. Consequently,
the corresponding Problem (MF-LQ) is not meaningful.
\eex

From the above, we see that before investigating Problem (MF-LQ), we
should find conditions for the system and the cost functional so
that the set $\cU_{ad}[0,\infty)$ is at least non-empty and
hopefully it admits an accessible characterization. To this end, let
us first look at the following uncontrolled linear MF-FSDE (which
amount to saying that taking $u(\cd)=0$ or letting $B=\bar B=D=\bar D=0$):
\bel{MF-state2}\left\{\2n\ba{ll}
\ns\ds dX(t)=\Big\{AX(t)+\bar A\dbE[X(t)]\Big\}dt+\Big\{CX(t)+\bar C\dbE[X(t)]\Big\}dW(t),\qq t\ge0, \\
\ns\ds X(0)=x,\ea\right.\ee
where $A,\bar A,C,\bar C\in\dbR^{n\times n}$ are given matrices. The
above uncontrolled system is briefly denoted by $[A,\bar A,C,\bar C]$.
For simplicity, we also denote $[A,C]=[A,0,C,0]$ (the linear
SDE without mean-fields), and $A=[A,0]\equiv[A,0,0,0]$ (the linear
ordinary differential equation, ODE, for short). Let us now
introduce the following definition.

\ms

\bde\label{def3.2}
\rm (i) System $[A,\bar A,C,\bar C]$ is said to
be $L^2$-{\it exponentially stable} if for any $x\in\dbR^n$, the
solution $X(\cd)\equiv X(\cd\,;x)\in\cX_{loc}[0,\infty)$ of
(\ref{MF-state2}) satisfies the following:
%
$$\lim_{t\to\infty}e^{\l t}\dbE|X(t)|^2=0,$$
for some $\l>0$.

\ms

(ii) System $[A,\bar A,C,\bar C]$ is said to be $L^2$-{\it globally
integrable} if for any $x\in\dbR^n$, the solution $X(\cd)\equiv
X(\cd\,;x)\in\cX_{loc}[0,\infty)$ of (\ref{MF-state2}) is in
$\cX[0,\infty)$, namely,
%
$$\int_0^\infty\dbE|X(t)|^2dt<\infty.$$

(iii) System $[A,\bar A,C,\bar C]$ is said to be $L^2$-{\it
asymptotically stable} if for any $x\in\dbR^n$, the solution
$X(\cd)\equiv X(\cd\,;x)\in\cX_{loc}[0,\infty)$ of (\ref{MF-state2})
satisfies the following:
\bel{lim=0}\lim_{t\to\infty}\dbE|X(t)|^2=0.\ee

(iv) Let (J) hold. System $[A,\bar A,C,\bar C]$ is said to be
$L^2_{Q,\bar Q}$-{\it globally integrable} if for any $x\in\dbR^n$,
the solution $X(\cd)\equiv X(\cd\,;x)\in\cX_{loc}[0,\infty)$ of
(\ref{MF-state2}) satisfies (\ref{QQ-stable}).
\ede

\ms

It is clear that the above (i)--(iii) can be defined for linear SDE
system $[A,C]=[A,0,C,0]$, and linear ODE system $A=[A,0]$. By a
standard result, we know that the above (i)--(iii) are equivalent
for linear ODEs. For general linear MF-SDEs, we have the following
result.

\bp\label{prop3.3}
\sl Among the following statements:

\ms

{\rm(i)} System $[A,\bar A,C,\bar C]$ is $L^2$-exponentially stable;

\ms

{\rm(ii)} System $[A,\bar A,C,\bar C]$ is $L^2$-globally integrable;

\ms

{\rm(iii)} System $[A,\bar A,C,\bar C]$ is $L^2$-asymptotically stable;

\ms

{\rm(iv)} Let {\rm(J)} hold. System $[A,\bar A,C,\bar C]$ is $L^2_{Q,\bar Q}$-globally integrable.

\ms

\no The following implications hold:
$$\ba{ll}
\ns\ds\hb{\rm(i)}~~\Ra~~\hb{\rm(ii)}~~\Ra~~\hb{\rm(iii)}; \\
\ns\ds\hb{when {\rm(J)} holds, \rm(ii)}~~\Ra~~\hb{\rm(iv)};\qq\hb{when {\rm(J)$'$} holds, \rm(iv)}~~\Ra~~{\rm(ii)}.\ea$$
\ep

\it Proof. \rm The implications (i) $\Ra$ (ii) and (ii) $\Ra$ (iv)
(under (J)) are clear. It is also clear that under (J)$'$, (iv)
$\Ra$ (ii). We now prove (ii) $\Ra$ (iii). Note that
$$\ba{ll}
\ns\ds\dbE|X(t)|^2 = \ds|x|^2+\dbE\int_0^t\(2\lan X(s),AX(s)+\bar A\dbE[X(s)]\ran+|CX(s)+\bar C\dbE[X(s)]|^2\)ds\\
\ns\ds\qq\qq \leq |x|^2+L\int_0^t\dbE|X(s)|^2ds\le|x|^2+L\dbE\int_0^\infty|X(s)|^2ds.\ea$$
Hereafter $L>0$ stands for a generic constant which could be
different from line to line. Thus, $\dbE|X(t)|^2$ is bounded
uniformly in $t\in[0,\infty)$. Consequently, for any
$0\le\t<t<\infty$,
$$\ba{ll}
\ns\ds\Big|\dbE|X(t)|^2-\dbE|X(\t)|^2\Big| \leq \ds\dbE\int_\t^t\(2|\lan X(s),AX(s)+\bar A\dbE[X(s)]\ran|+|CX(s)+\bar C\dbE[X(s)]|^2\)ds\\
\ns\ds\qq\qq\qq\qq\qq\1n \leq L(t-\tau).\ea$$
Hence, $t\mapsto\dbE|X(t)|^2$ is uniformly continuous on
$[0,\infty)$, which, together with the integrability of
$\dbE|X(\cd)|^2$ over $[0,\infty)$, leads to (\ref{lim=0}). \endpf

\ms

Let us make the following remarks.

\ms

$\bullet$ When (J) holds but (J)$'$ does not hold, the $L^2_{Q,\bar
Q}$-global integrability of the system does not imply the
$L^2$-global integrability of the system in general.

\ms

$\bullet$ It is not clear if (iii) implies (ii), although these two
are equivalent for ODE case.

\ms

$\bullet$ The notion that is the most relevant to our Problem
(MF-LQ) is the $L^2_{Q,\bar Q}$-global integrability.

\ms

Our next goal is to explore when (ii) implies (i). To this end, we
first look the case $\bar A=\bar C=0$. In this case, our system
becomes system $[A,C]$:
\bel{SDE}\left\{\2n\ba{ll}
\ns\ds dX(t)=AX(t)dt+CX(t)dW(t),\qq t\ge0,\\
\ns\ds X(0)=x.\ea\right.\ee
For such a system, instead of $L^2_{Q,\bar Q}$-global integrability,
we may introduce the following.

\ms

\bde\label{def3.4} \rm Let $Q\ge0$. System $[A,C]$ is said to be
$L^2_Q$-globally integrable if for any $x\in\dbR^n$, the solution
$X(\cd)\equiv X(\cd\,;x)$ of (\ref{SDE}) satisfies
$$\dbE\int_0^\infty\lan QX(t),X(t)\ran dt<\infty.$$
In the case that $Q>0$, the $L^2_Q$-global integrability is simply
called the $L^2$-global integrability which is equivalent to
$X(\cd\,;x)\in\cX[0,\infty)$ for all $x\in\dbR^n$.

\ede

We have the following result concerning the $L^2_Q$-global
integrability of $[A,C]$.

\ms

\bp\label{prop3.5} \sl Let $Q\ge0$. Then the following are
equivalent:

\ms

{\rm(i)} System $[A,C]$ is $L^2_Q$-globally integrable;

\ms

{\rm(ii)} The following Lyapunov equation admits a solution $P\ge0$:
\bel{Lyapunov0}PA+A^TP+C^TPC+Q=0;\ee

\ms

{\rm(iii)} The Lyapunov equation $(\ref{Lyapunov0})$ admits a
solution $P\in\cS^n$ and
$$\limsup_{t\to\infty}\dbE|X(t;x)|^2<\infty,\qq\forall x\in\dbR^n,$$
which is the case, in particular, if $[A,C]$ is $L^2$-asymptotically
stable.

\ms

In the above case, the solution $P$ of the above equation admits the
following representation:
\bel{P}
P=\dbE\int_0^\infty\F(t)^TQ\F(t)dt,\ee
where $\F(\cd)$ is the solution to the following:
$$\left\{\2n\ba{ll}
\ns\ds d\bar F(t)=A\bar F(t)dt+C\bar F(t)dW(t), \qq t\in[0,\infty), \\
\ns\ds\bar F(0)=I.\ea\right.$$ \ep

\rm

The above result should be standard. However, since the idea
contained in the proof will be useful below, for readers's
convenience, we present a proof here.

\ms

\it Proof. \rm (i) $\Ra$ (ii). Suppose system $[A,C]$ is
$L^2_Q$-globally integrable. We want to show that Lyapunov equation
(\ref{Lyapunov0}) admits a solution $P\ge0$. To this end, let us
consider the following linear ODE:
\bel{ODE}\left\{\2n\ba{ll}
\ns\ds-\dot\Th(t)+\Th(t)A+A^T\Th(t)+C^T\Th(t)C+Q=0,\qq t\in[0,\infty), \\
\ns\ds\Th(0)=0,\ea\right.\ee
which has a unique solution $\Th(\cd)$ defined on $[0,\infty)$. For
any fixed $\t>0$, we define
$$\bar\Th^\t(s)=\Th(\t-s),\qq s\in[0,\t].$$
Then $\bar\Th^\t(\cd)$ is the solution to the following:
$$\left\{\2n\ba{ll}
\ns\ds\dot{\bar\Th}^\tau(s)+\bar\Th^\t(s)A+A^T\bar\Th^\t(s)+C^T\bar\Th^\t(s)C+Q=0,\qq s\in[0,\t], \\
\ns\ds\bar\Th^\t(\t)=0.\ea\right.$$
For any $x\in\dbR^n$, let $X(\cd)\equiv X(\cd\,;x)$ be the solution
of (\ref{SDE}). Applying It\^o's formula to
$s\mapsto\lan\bar\Th^\t(s)X(s),X(s)\ran$, one has
$$\ba{ll}
\ns\ds-\lan\Th(\t)x,x\ran = \ds-\lan\bar\Th^\t(0)x,x\ran=\dbE\[\lan\bar\Th^\t(\t)X(\t),X(\t)\ran-\lan\bar\Th^\t(0)X(0),X(0)\ran\] \\
\ns\ds\qq\qq\qq = \dbE\int_0^\t\lan\big\{\dot{\bar\Th}^\t(s)+\bar\Th^\t(s)A+A^T\bar\Th^\t(s)+C^T\bar\Th^\t(s)C\big\}X(s),X(s)\ran ds \\
\ns\ds\qq\qq\qq = -\dbE\int_0^\t\lan QX(s),X(s)\ran ds=-\dbE\int_0^\t\lan\F(s)^TQ\F(s)x,x\ran ds.\ea$$
Thus, the solution $\Th(\cd)$ of (\ref{ODE}) admits the following representation:
$$\Th(\t)=\dbE\int_0^\t\F(s)^TQ\F(s)ds,\qq\t\ge0.$$
From the above, since $Q\ge0$, we see that $\t\mapsto\Th(\t)$ is
non-decreasing and by the $L^2_Q$-global integrability of $[A,C]$,
one has the following limit:
$$\lim_{\t\to\infty}\Th(\t)=\dbE\int_0^\infty\F(s)^TQ\F(s)ds\equiv P.$$
We claim that such a $P\ge0$ must be a solution to the Lyapunov
equation (\ref{Lyapunov0}). In fact, from (\ref{ODE}), one has
$$\Th(t)-\Th(t+1)+\(\int_t^{t+1}\Th(s)ds\)A+A^T\(\int_t^{t+1}\Th(s)ds\)+C^T\(\int_t^{t+1}\Th(s)ds\)C+Q=0.$$
Letting $t\to\infty$, we see that (\ref{Lyapunov0}) is satisfied by $P$.

\ms

(ii) $\Ra$ (i) Suppose there exists a $P\ge0$ satisfying (\ref{Lyapunov0}). Then
\bel{P0}\ba{ll}
\ns\ds\1n \dbE\lan PX(t),X(t)\ran-\lan Px,x\ran=\dbE\int_0^t\lan(PA+A^TP+C^TPC)X(s),X(s)\ran ds \\
\ns\ds\qq\qq\qq\qq\qq\qq=-\dbE\int_0^t\lan QX(s),X(s)\ran ds. \ea\ee
This implies
\bel{3.8}\dbE\int_0^t\lan QX(s),X(s)\ran ds
=\lan Px,x\ran-\dbE\lan PX(t),X(t)\ran\le\lan Px,x\ran,\qq t\ge0.\ee
Thus, the system is $L^2_Q$-globally integrable.

\ms

(i) $\Ra$ (iii) is clear.

\ms

(iii) $\Ra$ (i). Suppose (\ref{Lyapunov0}) has a solution
$P\in\cS^n$. Then by (\ref{3.8}), we have
$$\dbE\int_0^t\lan QX(s),X(s)\ran ds\le\lan Px,x\ran+L\dbE|X(t)|^2.$$
Hence, $[A,C]$ is $L^2_Q$-globally integrable. \endpf

\ms

Combining Propositions 3.3 and 3.5, we have the following result for
system $[A,C]$.

\ms

\bp\label{prop3.6}
\sl The following are equivalent:

\ms

{\rm(i)} System $[A,C]$ is $L^2$-exponentially stable;

\ms

{\rm(ii)} System $[A,C]$ is $L^2$-globally integrable;

\ms

{\rm(iii)} For any $Q>0$, the Lyapunov equation $(\ref{Lyapunov0})$
admits a solution $P>0$, and in this case, the representation
$(\ref{P})$ holds for this $P$;

\ms

{\rm(iv)} System $[A,C]$ is $L^2$-asymptotically stable, and for
some $Q>0$, Lyapunov equation $(\ref{Lyapunov0})$ admits a solution
$P\in\cS^n$.
\ep

\it Proof. \rm The implications (i) $\Ra$ (ii) $\Ra$ is clear. The
relations (ii) $\iff$ (iii) $\iff$ (iv) follow from Proposition
\ref{prop3.5}. The implication (iii) $\Ra$ (i) follows from
(\ref{P0}), together with the positive definiteness of $P$ and $Q$
and Gronwall's inequality.
\endpf

\ms

Now, let us return system $[A,\bar A,C,\bar C]$. We have the
following result.

\bt\label{thm3.7} \sl {\rm(i)} Suppose system $[A,\bar A,C,\bar C]$
is $L^2$-asymptotically stable. Then it is necessary that $A+\bar A$
is exponentially stable.

\ms

{\rm(ii)} If $A+\bar A$ is exponentially stable, then system
$[A,\bar A,C,\bar C]$ is $L^2$-exponentially stable if either
$[A,C]$ is $L^2$-globally integrable, or
\bel{C+C=0} C+\bar C=0. \ee
\et

\it Proof. \rm (i) Suppose (\ref{lim=0}) holds. Taking expectation
in (\ref{MF-state2}), we obtain
\bel{MF-state2E}\left\{\2n\ba{ll}
\ns\ds d\dbE[X(t)]=(A+\bar A)\dbE[X(t)]dt,\qq t\ge0, \\
\ns\ds\dbE[X(0)]=x.\ea\right.\ee
Hence,
$$\dbE[X(t)]=e^{(A+\bar A)t}x,\qq t\ge0.$$
Since
$$|\dbE[X(t)]|^2\le\dbE|X(t)|^2,\qq t\ge0,$$
the $L^2$-asymptotic stability of system $[A,\bar A,C,\bar C]$
implies the exponential stability of $A+\bar A$.

\ms

(ii) By (\ref{P(X-EX),X-EX0}) with $B=\bar B=D=\bar D=0$, we have,
for any $P\in\cS^n$,
\bel{2.35}\ba{ll}
\ns\ds\dbE\lan P\big\{X(t)-\dbE[X(t)]\big\},X(t)-\dbE[X(t)]\ran \\
\ns\ds\qq=\dbE\int_0^t\Big\{\lan(PA+A^TP+C^TPC)\big\{X(s)-\dbE[X(s)]\big\},X(s)-\dbE[X(s)]\ran \\
\ns\ds\qq\qq\qq+\lan P(C+\bar C)\dbE[X(s)],(C+\bar C)\dbE[X(s)]\ran\Big\}ds.\ea\ee
Hence, if (\ref{C+C=0}) holds, one has from the above that
$$\var[X(t)]\le L\int_0^t\var[X(s)]ds,\qq\forall t\ge0.$$
Then, by Gronwall's inequality, we obtain
$$\var[X(t)]=0,\qq t\ge0.$$
Consequently, if we let $2\l=-\max\si(A+\bar A)>0$, then
$$e^{2\l t}\dbE|X(t)|^2=e^{2\l t}\(\var[X(t)]+\big|\dbE[X(t)]\big|^2\)
=|e^{\l t}e^{(A+\bar A)t}x|^2\to0,\qq t\to\infty.$$
Thus, $[A,\bar A,C,\bar C]$ is $L^2$-exponentially stable.

\ms

Next, if $[A,C]$ is $L^2$-globally integrable, then by Proposition
\ref{prop3.6}, for $Q=I$, there exists a $P>0$ such that
%
$$PA+A^TP+C^TPC+I=0.$$
Hence, (\ref{2.35}) implies
$$\var[X(t)]\le-\m\int_0^t\var[X(s)]ds+L\int_0^t\big|\dbE[X(s)]\big|^2ds,\qq t\ge0,$$
for some $\m,L>0$, with $\m\ne\l=-\max\si(A+\bar A)>0$. By Gronwall's inequality,
$$\var[X(t)]\le L|x|^2\int_0^te^{-\m(t-s)}e^{-\l s}ds=L|x|^2{e^{-\l t}-e^{-\m t}\over\m-\l},\qq t\ge0.$$
This results in
$$\ba{ll}
\ns\ds\dbE|X(t)|^2=\var[X(t)]+\big|\dbE[X(t)]\big|^2\le L|x|^2{e^{-\l t}-e^{-\m t}\over\m-\l}+|e^{(A+\bar A)t}x|^2,\qq t\ge0.\ea$$
%
Therefore, the system $[A,\bar A,C,\bar C]$ is $L^2$-exponentially stable. This completes the proof.
\endpf

\ms

Note that the exponential stability of $A+\bar A$ together with the
$L^2$-global integrability of $[A,C]$ or (\ref{C+C=0}) are
sufficient conditions for the $L^2$-exponential stability of system
$[A,\bar A,C,\bar C]$. When $n=1$, these conditions are also
necessary in some sense. To be more precise, let us look at the
following one-dimensional system:
\bel{2.21}\left\{\2n\ba{ll}
\ns\ds dX(t)=\big\{aX(t)+\bar a\dbE[X(t)]\big\}dt+\big\{cX(t)+\bar c\dbE[X(t)]\big\}dW(t), \qq t\ge0, \\
\ns\ds X(0)=x.\ea\right.\ee
We have the following result.

\ms

\bp\label{prop3.8}
\sl For system $(\ref{2.21})$, the following are equivalent:

\ms

{\rm(i)} It is $L^2$-exponentially stable;

\ms

{\rm(ii)} It is $L^2$-globally integrable;

\ms

{\rm(iii)} It is $L^2$-asymptotically stable;

\ms

{\rm(iv)} $a+\bar a<0$, and
%
$$\hb{ either $2a+c^2<0$,~~~or~~~$2a+c^2\ge0$ and $c+\bar c=0$}.$$
\ep

\it Proof. \rm It suffices to prove the implication (iii)$\Ra$(iv).
By (\ref{PX,X0}) with $P=1$, $B=\bar B=D=\bar D=0$, $A=a$, $\bar
A=\bar a$, $C=c$, $\bar C=\bar c$, we have
$$\ba{ll}
\ns\ds\dbE|X(t)|^2 = x^2+\dbE\int_0^t\Big\{(2a+c^2)|X(s)|^2+(2\bar a+\bar c^2+2\bar cc)\big(\dbE[X(s)]\big)^2\Big\}ds \\
\ns\ds\qq\qq = x^2+\int_0^t\Big\{(2a+c^2)\dbE|X(s)|^2+[2\bar a-c^2+(c+\bar c)^2]x^2e^{2(a+\bar a)s}\Big\}ds.\ea$$
Thus,
$$\ba{ll}
\ns\ds\dbE|X(t)|^2 = e^{(2a+c^2)t}x^2+\big[2\bar a-c^2+(c+\bar c)^2\big]x^2\int_0^te^{(2a+c^2)(t-s)}e^{2(a+\bar a)s}ds \\
\ns\ds\qq\qq = e^{(2a+c^2)t}x^2+(2\bar a-c^2)x^2e^{(2a+c^2)t}\int_0^te^{(2\bar a-c^2)s}ds+(c+\bar c)^2x^2\int_0^te^{(2a+c^2)(t-s)}e^{2(a+\bar a)s}ds \\
\ns\ds\qq\qq = e^{(2a+c^2)t}x^2+x^2e^{(2a+c^2)t}\big[e^{(2\bar a-c^2)t}-1\big]+(c+\bar c)^2x^2\int_0^te^{(2a+c^2)(t-s)}e^{2(a+\bar a)s}ds \\
\ns\ds\qq\qq = x^2e^{2(a+\bar a)t}+(c+\bar c)^2x^2\int_0^te^{(2a+c^2)(t-s)}e^{2(a+\bar a)s}ds.\ea$$
Now, if (\ref{lim=0}) holds, then we must have
$$a+\bar a<0,$$
and
$$(c+\bar c)^2\int_0^te^{(2a+c^2)(t-s)}e^{2(a+\bar a)s}ds\to0.$$
Thus, under $a+\bar a<0$, if $c+\bar c\ne0$, then we need
$$\int_0^te^{(2a+c^2)(t-s)}e^{2(a+\bar a)s}ds=e^{(2a+c^2)t}\int_0^te^{(2\bar a-c^2)s}ds\to0.$$
Since $\int_0^te^{(2\bar a-c^2)s}ds$ is increasing, the above must
lead to $2a+c^2<0$. Also, if $2a+c^2\ge0$, we must have $c+\bar c=0$. This completes the proof.
\endpf

\ms

Now, for the $L^2_{Q,\bar Q}$-global integrability of system
$[A,\bar A,C,\bar C]$, we have the following result.

\ms

\bp\label{prop3.9}
\sl Let {\rm(J)} hold. If $[A,\bar A,C,\bar C]$
is $L^2_{Q,\bar Q}$-globally integrable, then $A+\bar A$ is
$L^2_{Q+\bar Q}$-globally integrable, i.e.,
\bel{2.41}\int_0^\infty|(Q+\bar Q)^{1\over2}e^{(A+\bar A)t}|^2dt<\infty.\ee
Conversely, if $(\ref{2.41})$ hold, then $[A,\bar A,C,\bar C]$ is
$L^2_{Q,\bar Q}$-globally integrable provided either (\ref{C+C=0})
holds, or $[A,C]$ is $L^2_Q$-globally integrable and
\bel{N(Q+Q)}\cN(Q+\bar Q)\subseteq\cN(C+\bar C),\ee
where $\cN(G)$ is the null space of $G$.
\ep

\it Proof. \rm Since,
$$\ba{ll}
\ns\ds\int_0^\infty\lan(Q+\bar Q)\dbE[X(t)],\dbE[X(t)]\ran dt\le\dbE\int_0^\infty\(\lan QX(t),X(t)\ran+\lan\bar Q\dbE[X(t)],\dbE[X(t)]\ran\)dt<\infty,\ea$$
we see that (\ref{2.41}) follows.

\ms

Next, let (\ref{2.41}) hold. If (\ref{C+C=0}) holds, we have (see
(\ref{2.35}) with $P=I$)
$$\ba{ll}
\ns\ds\2n \var[X(t)] = \dbE\lan\big\{X(t)-\dbE[X(t)]\big\},X(t)-\dbE[X(t)]\ran \\
\ns\ds\qq\qq = \dbE\int_0^t\Big\{\lan(A+A^T+C^TC)\big\{X(s)-\dbE[X(s)]\big\},X(s)-\dbE[X(s)]\ran+|(C+\bar C)\dbE[X(s)]|^2\Big\}ds \\
\ns\ds\qq\qq = \dbE\int_0^t\lan(A+A^T+C^TC)\big\{X(s)-\dbE[X(s)]\big\},X(s)-\dbE[X(s)]\ran ds\le L\int_0^t\var[X(s)]ds.\ea$$
Hence, by Gronwall's inequality, we obtain
$$\var[X(t)]=0,\qq t\ge0.$$
Consequently,
$$\ba{ll}
\ns\ds \dbE\int_0^\infty\(\lan QX(t),X(t)\ran+\lan\bar Q\dbE[X(t),\dbE[X(t)]\ran\)dt \\
\ns\ds = \dbE\int_0^\infty\(\lan Q\big\{X(t)-\dbE[X(t)\big\},X(t)-\dbE[X(t)]\ran+\lan(Q+\bar Q)\dbE[X(t)],\dbE[X(t)]\ran\)dt \\
\ns\ds \leq \int_0^\infty\(|Q|\var[X(t)]+|(Q+\bar Q)^{1\over2}e^{(A+\bar A)t}x|^2\)dt \\
\ns\ds = \int_0^\infty|(Q+\bar Q)^{1\over2}e^{(A+\bar A)t}x|^2dt<\infty,\ea$$
which gives the $L^2_{Q,\bar Q}$-global integrability.

\ms

Finally, if (\ref{2.41}) holds and $[A,C]$ is $L^2_Q$-globally
integrable, then by Proposition \ref{prop3.5}, we can find a $P\ge0$ solving
Lyapunov equation (\ref{Lyapunov0}). Let $X(\cd)$ be the solution of (\ref{MF-state2}).
Applying It\^o's formula to $\lan PX(\cd),X(\cd)\ran$, we get
%
$$\ba{ll}
\ns\ds \dbE\lan P\big\{X(t)-\dbE[X(t)]\big\},X(t)-\dbE[X(t)]\ran \\
\ns\ds = \dbE\int_0^t\Big\{\lan(PA+A^TP+C^TPC)\big\{X(s)-\dbE[X(s)]\big\},X(s)-\dbE[X(s)]\ran \\
\ns\ds\qq +\lan(C+\bar C)^TP(C+\bar C)\dbE[X(s)],\dbE[X(s)]\ran\Big\}ds \\
\ns\ds = \ds\dbE\2n\int_0^t\3n\Big\{\1n-\1n\lan Q\big\{X(s)\1n-\1n\dbE[X(s)]\big\},X(s)\1n-\1n\dbE[X(s)]\ran\1n+\1n\lan(C\2n+\1n\bar C)^T\1n P(C\2n+\1n\bar C)\dbE[X(s)],\dbE[X(s)]\ran\1n\Big\}ds.\ea$$
Now, condition (\ref{N(Q+Q)}) implies that
$$\lan P(C+\bar C)y,(C+\bar C)y\ran\le L\lan(Q+\bar Q)y,y\ran,\qq\forall y\in\dbR^n,$$
for some $L>0$. Thus,
$$\ba{ll}
\ns\ds \dbE\int_0^t\lan Q\big\{X(s)-\dbE[X(s)]\big\},X(s)-\dbE[X(s)]\ran ds \\
\ns\ds = \int_0^t\lan P(C+\bar C)\dbE[X(s)],(C+\bar C)\dbE[X(s)]\ran ds-\dbE\lan P\big\{X(t)-\dbE[X(t)]\big\},X(t)-\dbE[X(t)]\ran \\
\ns\ds \leq L\int_0^t\lan(Q+\bar Q)\dbE[X(s)],\dbE[X(s)]\ran ds.\ea$$
Consequently,
$$\ba{ll}
\ns\ds \dbE\int_0^\infty\(\lan QX(t),X(t)\ran+\lan\bar Q\dbE[X(t)],\dbE[X(t)]\ran\)dt\\
\ns\ds = \dbE\int_0^\infty\(\lan Q\big\{X(t)-\dbE[X(t)]\big\},X(t)-\dbE[X(t)\big\}\ran+\lan(Q+\bar Q)\dbE[X(t)],\dbE[X(t)]\ran\)dt \\
\ns\ds \leq (L+1)\int_0^\infty\lan(Q+\bar Q)\dbE[X(s)],\dbE[X(s)]\ran ds<\infty.\ea$$
This means that the system is $L^2_{Q,\bar Q}$-globally integrable.
\endpf

\ms

We point out that condition (\ref{N(Q+Q)}) holds if (\ref{C+C=0}) is
true or
%
$$Q+\bar Q>0.$$
Therefore, to have condition (\ref{N(Q+Q)}), we do not have to assume (J)$'$.

\section{MF-Stabilizability}

We now return to the controlled linear MF-FSDE (\ref{MF-state1})
which is denoted by $[A,\bar A,C,\bar C;B,\bar B,D,\bar D]$. With
this notation, we see that the uncontrolled MF-FSDE
(\ref{MF-state2}) is nothing but $[A,\bar A,C,\bar C;0,0,0,0]$. Note
also that in the case $\bar A=\bar C=0$ and $\bar B=\bar D=0$, the
system is a usual controlled linear SDE, which is simply denoted by
$[A,C;B,D]\equiv[A,0,C,0;B,0,D,0]$. Further, in the case $C=0$ and
$D=0$, the system is reduced to a classical controlled linear ODE,
which is denoted by $[A;B]\equiv[A,0,0,0;B,0,0,0]$. We now introduce
the following notion for general state equation (\ref{MF-state1}).

\ms

\bde\label{def4.1}
\rm (i) Let (J) hold. System $[A,\bar A,C,\bar C;B,\bar B,D,\bar D]$ is said to be MF-$L^2_{Q,\bar Q}$-{\it stabilizable} if there exists a pair $(K,\bar K)\in\dbR^{n\times m}\times\dbR^{n\times m}$ such that for any $x\in\dbR^n$ if
$X^{K,\bar K}(\cd)$ is the solution to the following:
%
$$\left\{\2n\ba{ll}
\ns\ds dX^{K,\bar K}(t) = \Big\{(A+BK)X^{K,\bar K}(t)+[\bar A+\bar B\bar K+B(\bar K-K)]\dbE[X^{K,\bar K}(t)]\Big\}dt\\
\ns\ds\qq\qq\qq +\Big\{(C+DK)X^{K,\bar K}(t)+[\bar C+\bar D\bar K+D(\bar K-K)]\dbE[X^{K,\bar K}(t)]\Big\}dW(t),\q t\ge0, \\
\ns\ds X^{K,\bar K}(0) = x,\ea\right.$$
and
\bel{u(K,K)}
u^{K,\bar K}(t)=K\big\{X^{K,\bar K}(t)-\dbE[X^{K,\bar K}(t)]\big\}+\bar K\dbE[X^{K,\bar K}(t)],\qq t\ge0,\ee
then
\bel{2.47}
\dbE\int_0^\infty\(\lan QX^{K,\bar K}(t),X^{K,\bar K}(t)\ran+ \lan\bar Q\dbE[X^{K,\bar K}(t)],\dbE[ X^{K,\bar K}(t)]\ran+|u^{K,\bar K}(t)|^2\)dt<\infty.\ee
In this case, the pair $(K,\bar K)$ is called an MF-$L^2_{Q,\bar
Q}$-{\it stabilizer} of the system. In the case that (\ref{2.47}) is
replaced by the following:
%
$$\dbE\int_0^\infty\(|X^{K,\bar K}(t)|^2+|u^{K,\bar K}(t)|^2\)dt<\infty,$$
we simply say that the system is MF-$L^2$-stabilizable, and $(K,\bar K)$ is called an MF-$L^2$-stabilizer of the system.

\ms

(ii) Let (J) hold. System $[A,\bar A,C,\bar C;B,\bar B,D,\bar D]$ is
said to be $L^2_{Q,\bar Q}$-{\it stabilizable} if there exists a
$K\in\dbR^{n\times m}$ such that for any $x\in\dbR^n$, if $X^K(\cd)$
is the solution to the following:
\bel{closed-loop0}\left\{\2n\ba{ll}
\ns\ds dX^K(t) = \Big\{(A+BK)X^K(t)+(\bar A+\bar BK)\dbE[X^K(t)]\Big\}dt \\
\ns\ds\qq\qq\q +\Big\{(C+DK)X^K(t)+(\bar C+\bar DK)\dbE[X^K(t)]\Big\}dW(t),\qq t\ge0, \\
\ns\ds X^K(0) = x,\ea\right.\ee
and
%
$$u^K(t)=KX^K(t),\qq t\ge0,$$
then
\bel{2.51}
\dbE\int_0^\infty\(\lan QX^K(t),X^K(t)\ran+\lan\bar Q\dbE[X^K(t)],\dbE[X^K(t)]\ran+|u^K(t)|^2\)dt<\infty.\ee
In this case, $K$ is called an $L^2_{Q,\bar Q}$-{\it stabilizer} of
the system. In the case that $\bar Q=0$, we simply say that the
system is $L^2_Q$-{\it stabilizable}, and $K$ is called an
$L^2_Q$-{\it stabilizer}. If (\ref{2.51}) is replaced by
%
$$\dbE\int_0^\infty\(|X^K(t)|^2+|u^K(t)|^2\)dt<\infty,$$
we further simply say that the system is $L^2$-{\it stabilizable},
and $K$ is called an $L^2$-{\it stabilizer} of the system.
\ede

The importance of the notions defined in the above definition is
that if (J) holds and $[A,\bar A,C,\bar C;B,\bar B,D,\bar D]$ is
MF-$L^2_{Q,\bar Q}$-stabilizable, then $\cU_{ad}[0,\infty)$ is
nonempty since $u^{K,\bar K}(\cd)$ defined by (\ref{u(K,K)}) is in
$\cU_{ad}[0,\infty)$. In particular, $\cU_{ad}[0,\infty)$ is
nonempty if the system $[A,\bar A,C,\bar C;B,\bar B,D,\bar D]$ is
MF-$L^2$-stabilizable.

\ms

It is seen that when system $[A,\bar A,C,\bar C;B,\bar B,D,\bar D]$
is MF-$L^2_{Q,\bar Q}$-stabilizable, then the uncontrolled system
$[A+BK,\bar A+\bar B\bar K+B(\bar K-K),C+DK,\bar C+\bar D\bar
K+D(\bar K-K)]$ is $L^2_{Q,\bar Q}$-globally integrable. Also,
system $[A,\bar A,C,\bar C;B,\bar B,D,\bar D]$ is
$L^2_Q$-stabilizable if and only if
$$\dbE\int_0^\infty\(\lan QX^K(t),X^K(t)\ran+|u^K(t)|^2\)dt<\infty.$$
Moreover, it is clear that the $L^2$-stabilizability of system
$[A,C;B,D]$ we defined here is the classic stabilizability of the
controlled SDE system.

\ms

Note that system (\ref{MF-state1}) is $L^2_{Q,\bar Q}$-stabilizable
(resp. $L^2$-stabilizability) if it is MF-$L^2_{Q,\bar Q}$-stabilizable (resp. MF-$L^2$-stabilizability) with $K=\bar K$.
Therefore, the former is a special case of the later. The following
example shows that in general, the MF-$L^2$-stabilizability does not
imply the $L^2$-stabilizability.

\ms

\bex\label{ex4.2}
\rm Consider the following one-dimensional controlled MF-FSDE:
%
$$\left\{\2n\ba{ll}
\ns\ds dX(t) = \big\{aX(t)+\bar a\dbE[X(t)]+bu(t)+\bar b\dbE[u(t)]\big\}dt \\
\ns\ds\qq\qq +\big\{cX(t)+\bar c\dbE[X(t)]+du(t)+\bar d\dbE[u(t)]\big\}dW(t),\qq t\ge0, \\
\ns\ds X(0) = x.\ea\right.$$
Suppose the above system is MF-$L^2$-stabilizable. Then, there are
$k,\bar k\in\dbR$ such that with
$$u(t)=k\big\{X(t)-\dbE[X(t)]\big\}+\bar k\dbE[X(t)],\qq t\ge0,$$
the closed-loop system:
$$\ba{ll}
\ns\ds dX(t) = \big\{(a+bk)X(t)+[\bar a+\bar b\bar k+b(\bar k-k)]\dbE[X(t)]\big\}dt \\
\ns\ds\qq\qq +\big\{(c+dk)X(t)+[\bar c+\bar d\,\bar k+d(\bar k-k)]\dbE[X(t)]\big\}dW(t),\qq t\ge0,\ea$$
is $L^2$-globally integrable. By Proposition \ref{prop3.8}, this is equivalent to the following:
$$a+\bar a+(b+\bar b)\bar k<0,$$
and either
$$2(a+bk)+(c+dk)^2<0,$$
or
$$2(a+bk)+(c+dk)^2\ge0,\q c+\bar c+(d+\bar d)\bar k=0.$$
Let
$$b+\bar b=1,\q d=1,\q\bar d=-1,\q c+\bar c\ne0.$$
Then we need and only need
\bel{2.33}a+\bar a+\bar k\equiv-\l<0,\qq2(a+bk)+(c+k)^2<0,\ee
for some $k,\bar k\in\dbR$. The first condition in (\ref{2.33}) can
always be achieved. The second one is equivalent to the following:
$$0>k^2+2(b+c)k+2a+c^2=(k+b+c)^2+2a+c^2-(b+c)^2,$$
which is possible if
\bel{2.28}2a+c^2-(b+c)^2<0.\ee
On the other hand, in order the system to be stabilizable, we need
$k=\bar k$, and
$$\ba{ll}
\ns\ds 0 > \ds2(a+b\bar k)+(c+\bar k)^2=2\big[a-b(a+\bar a+\l)\big]+\big[c-(a+\bar a+\l)\big]^2 \\
\ns\ds\q\1n = \ds\l^2-2(a+\bar a+b-c)\l+(a+\bar a-c)^2+2[a-b(a+\bar a)],\ea$$
for some $\l>0$. This is impossible if, say,
\bel{2.29}
c\ge a+\bar a+b,\qq a-b(a+\bar a)\ge 0.\ee
It is easy to find cases that (\ref{2.28})--(\ref{2.29}) hold.
Hence, we see that MF-$L^2$-stabilizability does not imply
$L^2$-stabilizability, in general.
\eex

Now, we present a result concerning the MF-$L^2_{Q,\bar Q}$-stabilizability of system (\ref{MF-state1}).

\ms

\bt\label{thm4.3}
\sl Let {\rm(J)} hold.

\ms

{\rm(i)} If system $(\ref{MF-state1})$ is MF-$L^2_{Q,\bar Q}$-stabilizable,
then the controlled ODE system $[A+\bar A;B+\bar B]$ is $L^2_{Q+\bar Q}$-stabilizable,
i.e., for some $\bar K\in\dbR^{m\times n}$,
\bel{2.57}\int_0^\infty|(Q+\bar Q)^{1\over2}e^{[A+\bar A+(B+\bar B)\bar K]t}|^2dt<\infty.\ee

\ms

{\rm(ii)} Suppose the following holds for some $\bar K\in\dbR^{m\times n}$ satisfying $(\ref{2.57})$:
\bel{2.58}
\cN(Q+\bar Q)\subseteq\cN\big(C+\bar C)+(D+\bar D)\bar K\big).\ee
Further, suppose the controlled SDE system $[A,C;B,D]$ is
$L^2_Q$-stabilizable. Then the controlled MF-FSDE system $[A,\bar A,C,\bar C;B,\bar B,D,\bar D]$ is MF-$L^2_{Q,\bar Q}$-stabilizable.

\ms

{\rm(iii)} Suppose the following holds for some $\bar
K\in\dbR^{n\times m}$ satisfying $(\ref{2.57})$:
\bel{C+C=00}
C+\bar C+(D+\bar D)\bar K=0.\ee
Then the controlled MF-FSDE system $[A,\bar A,C,\bar C;B,\bar
B,D,\bar D]$ is MF-$L^2_{Q,\bar Q}$-stabilizable.
\et

\it Proof. \rm Under (\ref{u(K,K)}), the closed-loop system takes
form (\ref{closed-loop0}). According to Proposition \ref{prop3.9}, we know that
if (\ref{closed-loop0}) is $L^2_{Q,\bar Q}$-globally integrable, it
is necessary that (\ref{2.57}) holds, which proves (i). Further,
when (\ref{2.57}) holds, the system (\ref{closed-loop0}) is
$L^2_{Q,\bar Q}$-globally integrable if either the system
$[A+BK,C+DK]$ is stable and (\ref{2.58}) holds, which proves (ii),
or (\ref{C+C=00}) holds with the same $\bar K$ which proves (iii).
\endpf

\ms

The above leads to the following corollary.

\ms

\bc\label{cor4.4}
\sl {\rm(i)} If system $(\ref{MF-state1})$ is
MF-$L^2$-stabilizable, then the controlled ODE system $[A+\bar
A;B+\bar B]$ is stabilizable, i.e., there exists a $\bar
K\in\dbR^{n\times m}$ such that
\bel{A+A+(B+B)K<0}
\si\big(A+\bar A+(B+\bar B)\bar K\big)\subseteq\dbC^-.\ee

\ms

{\rm(ii)} Suppose controlled ODE system $[A+\bar A;B+\bar B]$ is
stabilizable, and controlled SDE system $[A,C;B,D]$ is
$L^2$-stabilizable. Then the controlled MF-FSDE system
$[A,\bar A,C,\bar C;B,\bar B,D,\bar D]$ is MF-$L^2$-stabilizable.

\ms

{\rm(iii)} Suppose $(\ref{C+C=00})$ holds for some $\bar
K\in\dbR^{n\times m}$ satisfying $(\ref{A+A+(B+B)K<0})$. Then the
controlled MF-FSDE system $[A,\bar A,C,\bar C;B,\bar B,D,\bar D]$ is
MF-$L^2$-stabilizable.
\ec

\rm

Note that conditions assumed in (ii) of Corollary \ref{cor4.4} do not
involve $\bar C$ and $\bar D$. However, condition (\ref{C+C=00})
involves both $\bar C$ and $\bar D$. We point out that (\ref{C+C=00}) means that
\bel{C in D}
\cR(C+\bar C)\subseteq\cR(D+\bar D).\ee
In the case that $m<n$, the above could be a big restriction on
$C+\bar C$ and $D+\bar D$. Moreover, we have to find the same $\bar
K\in\dbR^{m\times n}$ such that (\ref{A+A+(B+B)K<0}) and
(\ref{C+C=00}) hold at the same time. If we let $(D+\bar D)^+$ be
the Moore-Penrose pseudo-inverse of $D+\bar D$ (\cite{Ben-Israel 2003}), then the solution of (\ref{C+C=00}) is given by
$$\bar K=-(D+\bar D)^+(C+\bar C)+\big[I-(D+\bar D)^+(D+\bar D)\big]\wt K,$$
for some $\wt K\in\dbR^{m\times n}$. Thus, we need
$$\si\(A+\bar A+(B+\bar B)\big\{-(D+\bar D)^+(C+\bar C)+\big[I-(D+\bar D)^+(D+\bar D)\big]\wt K\big\}\)\subseteq\dbC^-,$$
for some $\wt K\in\dbR^{m\times n}$, which means the ODE system
\bel{2.20}
\Big[A+\bar A-(B+\bar B)(D+\bar D)^+(C+\bar C);(B+\bar B)\big[I-(D+\bar D)^+(D+\bar D)\big]\Big]\ee
is stabilizable. Hence, we obtain the following result.

\ms

\bp\label{prop4.5}
\sl Let $(\ref{C in D})$ hold. Then $[A,\bar
A,C,\bar C;B,\bar B,D,\bar D]$ is MF-$L^2$-stabilizable if ODE
system $(\ref{2.20})$ is stabilizable, which is the case, if, in
particular, $m=n$, $D+\bar D$ is invertible, and
\bel{2.42}\si\big(A+\bar A-(B+\bar B)(D+\bar D)^{-1}(C+\bar C)\big)\subseteq\dbC^-.\ee
\ep

\rm

Condition (\ref{2.42}) seems that the MF-$L^2$-stabilizability of
system $[A,\bar A,C,\bar C;B,\bar B,D,\bar D]$ could be nothing to
do with the stabilizability of the controlled linear SDE system
$[A,C;B,D]$. However, in the case that $\bar A=\bar C=0$ and $\bar
B=\bar D=0$, we have the following controlled linear SDE:
%
$$dX(t)=\Big\{AX(t)+Bu(t)\Big\}dt+\Big\{CX(t)+Du(t)\Big\}dW(t),\qq t\ge0.$$
Suppose $m=n$ and $D^{-1}$ exists. Then condition (\ref{2.42})
becomes
\bel{2.36} \si\big(A-BD^{-1}C\big)\subseteq\dbC^-.\ee
In this case, if we take
$$u(t)=-D^{-1}CX(t),\qq t\ge0,$$
then the closed-loop system becomes
$$dX(t)=(A-BD^{-1}C)X(t)dt,\qq t\ge0,$$
which is stable if (\ref{2.36}) holds. Interestingly, if we let
\bel{2.37}
\bar u(t)=-D^{-1}C\dbE[X(t)],\qq\forall t\ge0,\ee
then the closed-loop system reads
$$\left\{\2n\ba{ll}
\ns\ds dX(t) = \Big\{AX(t)-BD^{-1}C\dbE[X(t)]\Big\}dt+C\Big\{X(t)-\dbE[X(t)]\Big\}dW(t),\qq t\ge0, \\
\ns\ds X(0) = x.\ea\right.$$
It is not hard to see that the unique solution $X(\cd)$ of the above
is deterministic and given by
$$X(t)=e^{(A-BD^{-1}C)t}x,\qq t\ge0.$$
Therefore the system is also asymptotically stable under feedback
control (\ref{2.37}).

\ms

\section{Stochastic LQ Problems}

In this section, we study a classic stochastic LQ problem, which
will be crucial for Problem (MF-LQ). We consider the following controlled SDE:
$$\left\{\2n\ba{ll}
\ns\ds dX(t) = \Big\{AX(t)+Bu(t)\Big\}dt+\Big\{CX(t)+Du(t)\Big\}dW(t),\qq t\ge0, \\
\ns\ds X(0) = x,\ea\right.$$
and cost functional
$$J^0(x;u(\cd))=\dbE\int_0^\infty\big\{\lan QX(t),X(t)\ran+\lan Ru(t),u(t)\ran\big\}dt.$$
Let
$$\left\{\2n\ba{ll}
\ds\cX_{ad}^Q[0,\infty)=\Big\{X(\cd)\in\cX_{loc}[0,\infty)\Big|\dbE\int_0^\infty\lan QX(t),X(t)\ran dt<\infty\Big\}, \\
\ns\ds\cU_{ad}^Q[0,\infty)=\Big\{u(\cd)\in\cU[0,\infty)\Big|X(x;u(\cd))\in\cX_{ad}^Q[0,\infty),\q \forall x\in\dbR^n\Big\}.\ea\right.$$

\subsection{A Classic Stochastic LQ Problem}

We introduce the following assumptions.

\ms

{\bf(J)$^*$} The matrices $Q\in\cS^n$ and $R\in\cS^m$ satisfy
$$Q\ge0,\qq R>0.$$

\ms

{\bf(S)$^*$} The system $[A,C;B,D]$ is $L^2_Q$-stabilizable.

\ms

Let us pose the following problem.

\ms

\bf Problem (LQ). \rm For any $x\in\dbR^n$, find a
$u_*(\cd)\in\cU_{ad}^Q[0,\infty)$ such that
$$J^0(x;u_*(\cd))=\inf_{u(\cd)\in\cU_{ad}^Q[0,\infty)}J^0(x;u(\cd))=V^0(x).$$

\ms

We have the following result.

\ms

\bt\label{thm5.1} \sl Let {\rm(J)$^*$} and {\rm(S)$^*$} hold. Then
Problem (LQ) admits a unique optimal control
$u_Q(\cd)\in\cU_{ad}^Q[0,\infty)$. Moreover, the following ARE
admits a solution $P\ge0$:
$$PA+A^TP+C^TPC+Q-(PB+C^TPD)(R+D^TPD)^{-1}(B^TP+D^TPC)=0,$$
and $\G$ is an $L^2_Q$-stabilizer of $[A,C;B,D]$, where
\bel{Gamma}
\G=-(R+D^TPD)^{-1}(B^TP+D^TPC).\ee
Further, the optimal control $u_Q(\cd)$ is given by
$$u^Q(t)=\G X^Q(t),\qq t\ge0,$$
with the optimal state process $X^Q(\cd)$ being the solution of
closed-loop system:
$$\left\{\2n\ba{ll}
\ns\ds dX^Q(t) = (A+B\G)X^Q(t)dt+(C+D\G)X^Q(t)dW(t),\qq t\ge0, \\
\ns\ds X^Q(0) = x,\ea\right.$$
and
\bel{V=Pxx}
\lan Px,x\ran=\inf_{u(\cd)\in\cU_{ad}^Q[0,\infty)}J^0(x;u(\cd))\equiv V^0(x),\qq\forall x\in\dbR^n.\ee
\et

\it Proof. \rm First of all, it is clear that under (J)$^*$ and
(S)$^*$, the set $\cU_{ad}^Q[0,\infty)$ is nonempty, and
$(x,u(\cd))\mapsto J^0(x;u(\cd))$ is a quadratic functional,
coercive with respect to $u(\cd)\in\cU_{ad}^Q[0,\infty)$. Thus for
any $x\in\dbR^n$, there exists a unique optimal control
$u^Q(\cd)\in\cU_{ad}^Q[0,\infty)$, and the value function $x\mapsto
V^0(x)$ must be of form (\ref{V=Pxx}) for some $P\ge0$. We now would
like to determine $P$ and the optimal pair $(X_*(\cd),u_*(\cd))$. To
this end, let us introduce
$$J^0_T(x;u(\cd))=\dbE\int_0^T\big\{\lan QX(t),X(t)\ran+\lan Ru(t),u(t)\ran\big\}dt,\qq T>0,$$
where $u(\cd)\in\cU_{loc}[0,\infty)$ and $X(\cd)=X(\cd\,;x,u(\cd))$.
It is standard that under (J)$^*$, there exists a unique
$u^Q_T(\cd)\in\cU[0,T]$ such that
$$V^0_T(x)\equiv\inf_{u(\cd)\in\cU[0,T]}J^0_T(x;u(\cd))=J^0_T(x;u^Q_T(\cd))=\lan P(0;T)x,x\ran,\qq\forall x\in\dbR^n,$$
with $P(\cd\,;T)$ being the solution to the following differential
Riccati equation:
\bel{Riccati-T0}\left\{\2n\ba{ll}
\ns\ds\dot P(t;T)+P(t;T)A+A^TP(t;T)+C^TP(t;T)C+Q \\
\ns\ds\q-\big[P(t;T)B\2n+\1n C^T\1n P(t;T)D\big] \big[R\2n+\1n D^T\1n P(t;T)D\big]^{-1}\big[B^T\1n P(t;T)\1n+\1n D^T\1n P(t;T)C\big]\1n=\1n0,\q t\in[0,T], \\
\ns\ds P(T;T)=0.\ea\right.\ee
Moreover, the optimal control $u_T(\cd)$ can be represented as
follows:
%
$$u^Q_T(t)=\G(t;T)X^Q_T(t),\qq t\in[0,T],$$
with
%
$$\G(t;T)=-\big[R+D^TP(t;T)D\big]^{-1}\big[B^TP(t;T)+D^TP(t;T)C\big],\qq t\in[0,T],$$
and $X^Q_T(\cd)$ is the solution to the following closed-loop system:
\bel{closed-T}\left\{\2n\ba{rcl}
\ns\ds dX^Q_T(t) & \5n=\5n & \big[A+B\G(t;T)\big]X^Q_T(t)dt+\big[C+D\G(t;T)\big]X^Q_T(t)dW(t),\qq t\in[0,T], \\
\ns\ds X^Q_T(0) & \5n=\5n & x.\ea\right.\ee
%
Now, it is clear that
$$J^0_T(x;u(\cd))\le J^0_{\bar T}(x;u(\cd)),\qq\forall u(\cd)\in\cU[0,\bar T],\qq 0\le T\le\bar T<\infty.$$
Therefore, one has
%
$$0\le P(0;T)\le P(0;\bar T),\qq\forall0\le T\le\bar T<\infty.$$
On the other hand, since
$$\cU^Q_{ad}[0,T]\equiv\Big\{u(\cd)\big|_{[0,T]}\Bigm|u(\cd)\in\cU^Q_{ad}[0,\infty)\Big\}\subseteq\cU[0,T],$$
it is true that
$$\ba{ll}
\ns\ds\lan P(0;T)x,x\ran\equiv V^0(x)=\inf_{u(\cd)\in\cU[0,T]}J^0_T(x;u(\cd))\le\inf_{u(\cd)\in\cU^Q_{ad}[0,T]}J^0_T(x;u(\cd)) \\
\ns\ds\qq\qq\qq\le\inf_{u(\cd)\in\cU^Q_{ad}[0,\infty)}J^0(x;u(\cd))=V^0(x)\equiv\lan Px,x\ran,\qq\forall x\in\dbR^n.\ea$$
Combining the above, we see that
$$0\le P(0;T)\le P(0;\bar T)\le P,\qq\forall0\le T\le\bar T<\infty.$$
This implies that
\bel{lim P(T)}
\lim_{T\to\infty}P(0;T)=\bar P\le P,\ee
for some $\bar P(\cdot)\ge0$. Now, we introduce the following differential
Riccati equation (on $[0,\infty)$):
$$\left\{\2n\ba{ll}
\ns\ds-\dot{\bar P}(s)+\bar P(s)A+A^T\bar P(s)+C^T\bar P(s)C+Q \\
\ns\ds\q-\big[\bar P(s)B+C^T\bar P(s)D\big]\big[R+D^T\bar P(s)D\big]^{-1}\big[B^T\bar P(s)+D^T\bar P(s)C\big]=0,\q s\ge0, \\
\ns\ds\bar P(0)=0.\ea\right.$$
For any $T>0$, let
$$\wt P(t;T)=\bar P(T-t),\qq t\in[0,T].$$
Then by the uniqueness, we must have
$$P(t;T)=\wt P(t;T)=\bar P(T-t),\qq\qq t\in[0,T].$$
Hence,
$$P(0;T)=\bar P(T),\qq\qq T\ge0.$$
From (\ref{lim P(T)}), we have
\bel{lim Th}
\lim_{t\to\infty}\bar P(t)=\bar P.\ee
This $\bar P\ge0$ must be a solution to the algebraic Riccati equation:
$$\bar PA+A^T\bar P+C^T\bar PC-(\bar PB+C^T\bar PD)(R+D^T\bar PD)^{-1}(B^T\bar P+D^T\bar PC)+Q=0.$$
Further, from (\ref{lim Th}), one has
$$\lim_{T\to\infty}P(t;T)=\lim_{T\to\infty}\bar P(T-t)=\bar P,\qq t\ge0.$$
Consequently,
\bel{Gamma0}
\lim_{T\to\infty}\G(t;T)=-(R+D^T\bar PD)^{-1}(B^T\bar P+D^T\bar PC)\equiv\G_0,\qq\forall t\ge0.\ee
Note that (suppressing $(t;T)$)
$$\ba{ll}
\ns\ds P(A+B\G)+(A+B\G)^TP+(C+D\G)^TP(C+D\G)+\G^TR\G \\
\ns\ds = P\big[A-B(R+D^TPD)^{-1}(B^TP+D^TPC)\big]+\big[A-B(R+D^TPD)^{-1}(B^TP+D^TPC)\big]^TP \\
\ns\ds\q +\big[C-D(R+D^TPD)^{-1}(B^TP+D^TPC)\big]^TP\big[C-D(R+D^TPD)^{-1}(B^TP+D^TPC)\big] \\
\ns\ds\q +(PB+CPD^T)(R+D^TPD)^{-1}R(R+D^TPD)^{-1}(B^TP+D^TPC) \\
\ns\ds = PA+A^TP+C^TPC-PB(R+D^TPD)^{-1}(B^TP+D^TPC) \\
\ns\ds\q -(PB+C^TPD)(R+D^TPD)^{-1}B^TP-(PB+C^TPD)(R+D^TPD)^{-1}D^TPC \\
\ns\ds\q -C^TPD(R+D^TD)^{-1}(B^TP+D^TPC) \\
\ns\ds\q +(PB+C^TPD)(R+D^TPD)^{-1}D^TPD(R+D^TPD)^{-1}(B^TP+D^TPC) \\
\ns\ds\q +(PB+CPD^T)(R+D^TPD)^{-1}R(R+D^TPD)^{-1}(B^TP+D^TPC) \\
\ns\ds = PA+A^TP+C^TPC-(PB+C^TPD)(R+D^TPD)^{-1}(B^TP+D^TPC) \\
\ns\ds\q -(PB+C^TPD)(R+D^TPD)^{-1}(B^TP+D^TPC) \\
\ns\ds\q +(PB+CPD^T)(R+D^TPD)^{-1}(B^TP+D^TPC) \\
\ns\ds = PA+A^TP+C^TPC-(PB+C^TPD)(R+D^TPD)^{-1}(B^TP+D^TPC).\ea$$
Next, we rewrite the differential Riccati equation
(\ref{Riccati-T0}) as follows:
$$\left\{\2n\ba{ll}
\ns\ds\dot P(t;T)+P(t;T)\big[A+B\G(t;T)\big]+\big[A+B\G(t;T)\big]^T\1n P(t;T) \\
\ns\ds\q+\big[C\2n+\1n D\G(t;T)\big]^T\1n P(t;T)\big[C\2n+\1n D\G(t;T)\big]\1n+\1n\G(t;T)^T\1n R\,\G(t;T)\1n+Q\1n=0,\qq t\in[0,T], \\
\ns\ds P(T;T)=0.\ea\right.$$
It is clear that (see (\ref{closed-T}) and (\ref{Gamma0}))
$$\lim_{T\to\infty}X_T^Q(t)=\bar X^Q(t),\qq t\ge0,$$
with $\bar X^Q(\cd)$ being the solution to the following:
$$\left\{\2n\ba{ll}
\ns\ds d\bar X^Q(t)=(A+B\G_0)\bar X^Q(t)dt+(C+D\G_0)\bar X^Q(t)dW(t),\qq t\ge0, \\
\ns\ds\bar X(0)=0.\ea\right.$$
Further,
$$\ba{ll}
\ns\ds\1n\lan P(0;T)x,x\ran = -\dbE\2n\int_0^T\3n\Big\{\1n\lan\1n\big\{\1n\dot P(t;T)\1n+\1n P(t;T)\big[A\1n+\1n B\G(t;T)\big]\1n+\1n\big[A(t;T)\1n+\1n B\G(t;T)\big]^TP(t;T) \\
\ns\ds\qq\qq\qq\q +\big[C\1n+\1n D\G(t;T)\big]^TP(t;T)\big[C\1n+\1n D\G(t;T)\big]\big\}X^Q_T(t),X^Q_T(t)\ran\1n\Big\}dt \\
\ns\ds\qq\qq\qq = \dbE\int_0^T\lan\big[Q+\G(t;T)^TR\G(t;T)\big]X^Q_T(t),X^Q_T(t)\ran dt \\
\ns\ds\qq\qq\qq = \dbE\int_0^T\(\lan QX^Q_T(t),X^Q_T(t)\ran+\lan R\G(t;T)X^Q_T(t),\G(t;T)X^Q_T(t)\ran\)dt.\ea$$
Thus, by Fatou's Lemma, we obtain (see also (\ref{lim P(T)}))
$$\lan Px,x\ran\ge\lan\bar Px,x\ran\ge\dbE\int_0^\infty\(\lan Q\bar X^Q(t),\bar X^Q(t)\ran+\lan R\G_0\bar X^Q(t),\G_0\bar X^Q(t)\ran\)dt\ge V^0(x)=\lan Px,x\ran,$$
which implies
$$\bar P=P,\qq\G_0=\G,$$
and $\bar X^Q(\cd)\in\cX^Q_{ad}[0,\infty)$. Also, $\G$ defined by
(\ref{Gamma}) is an $L^2_Q$-stabilizer of $[A,C;B,D]$. This
completes the proof. \endpf

\subsection{Stochastic MF-LQ Problem}

Having the above, let us now return to Problem (MF-LQ). We introduce
the following assumption.

\ms

{\bf(S)} The controlled ODE system $[A+\bar A;B+\bar B]$ is
stabilizable, and the controlled SDE system $[A,C;B,D]$ is
$L^2$-stabilizable.

\ms

From Corollary \ref{cor4.4}, we know that under (J) and (S), the system
$[A,\bar A,C,\bar C;B,\bar B,D,\bar D]$ is MF-$L^2$-stabilizable. We point out that it is possible for us to relax
(S) in various ways. However, for the simplicity of presentation, we
would like to keep the above (S). Let us first present the following
result.

\ms

Now, we are ready to state and prove the main result of this paper.

\ms

\bt\label{thm5.2} \sl
Let {\rm(J)} and {\rm(S)} hold. Then Problem
(MF-LQ) admits a unique optimal control
$u_*(\cd)\in\cU_{ad}[0,\infty)$, and the following AREs:
\bel{A-Riccati}\left\{\2n\ba{ll}
\ns\ds PA\1n+\1n A^T\1n P\1n+\1n C^T\1n PC\2n+\1n Q\1n-\1n(PB\2n+\1n C^TPD)(R+D^TPD)^{-1}(B^T\1n P\2n+\1n D^T\1n PC)=0, \\
\ns\ds\Pi(A\1n+\1n\bar A)\1n+\1n(A\1n+\1n\bar A)^T\Pi\1n+\1n(C\1n+\1n\bar C)^T\1n P(C\1n+\1n\bar C)\1n+\1n Q\1n+\1n\bar Q \\
\ns\ds\q-\big[\Pi(B\1n+\1n\bar B)\1n+\1n(C\1n+\1n\bar C)^T\1n P (D\1n+\1n\bar D)\big]\big[R+\bar R+(D+\bar D)^TP(D+\bar D)\big]^{-1} \\
\ns\ds\q\cd\big[(B\1n+\1n\bar B)^T\Pi\1n+\1n(D\1n+\1n\bar D)^T\1n P(C\1n+\1n\bar C)\big]=0,\ea\right.\ee
admits a solution pair $(P,\Pi)\in\bar\cS^n_+\times\bar\cS^n_+$.
Define
%
$$\left\{\2n\ba{ll}
\ns\ds\G=-(R+D^TPD)^{-1}(B^TP+D^TPC), \\
\ns\ds\bar\G=-\big[R+\bar R+(D+\bar D)^TP(D+\bar D)\big]^{-1}\big[(B+\bar B)^T\Pi +(D+\bar D)^TP(C+\bar C)\big].\ea\right.$$
Then $(\G,\bar\G)$ is an MF-$L^2_{Q,\bar Q}$-stabilizer of the
system. If $X_*(\cd)$ is the solution to the following MF-FSDE:
%
$$\left\{\2n\ba{ll}
\ns\ds dX_*(t) = \Big\{(A+B\G)X_*(t)+\big[\bar A+\bar B\bar\G+B(\bar\G-\G)\big]\dbE[X_*(t)]\Big\}dt \\
\ns\ds\qq\qq +\Big\{(C+D\G)X_*(t)+\big[\bar C+\bar D\bar\G+D(\bar \G-\G)\big]\dbE[X_*(t)]\Big\}dW(t),\qq t\ge0, \\
\ns\ds X_*(0) = x,\ea\right.$$
then
\bel{J=Pi}\ba{ll}
\ns\ds\inf_{u(\cd)\in\cU_{ad}[0,\infty)}J(x;u(\cd))=J(x;u_*(\cd))=\lan\Pi
x,x\ran,\qq\forall x\in\dbR^n,\ea\ee
with the optimal control $u_*(\cd)\in\cU_{ad}[0,\infty)$ admits the
following state feedback representation:
%
$$u_*(t)=\G\big\{X_*(t)-\dbE[X_*(t)]\big\}+\bar\G\dbE[X_*(t)],\qq t\ge0.$$
\et

\it Proof. \rm We know that under (J) and (S), the set
$\cU_{ad}[0,\infty)$ is nonempty, and convex. For any
$(x,u(\cd))\in\dbR^n\times\cU_{ad}[0,\infty)$, let
$X(\cd)=X(\cd\,;x,u(\cd))\in\cX[0,\infty)$. Then $J(x;u(\cd))$ is
well-defined and
%
$$\ba{ll}
\ns\ds\2n J(x;u(\cd)) = \dbE\2n\int_0^\infty\3n\Big\{\1n\lan QX(t),X(t)\ran\1n+\1n\lan\bar Q\dbE[X(t)],\dbE[X(t)]\ran\1n+\1n\lan Ru(t),u(t)\ran\1n+\1n\lan\bar R\dbE[u(t)],\dbE[u(t)]\ran\Big\}dt \\
\ns\ds\qq\qq = \ds\dbE\int_0^\infty\Big\{\lan Q\big\{X(t)-\dbE[X(t)]\big\},X(t)-\dbE[X(t)]\ran+\lan(Q+\bar Q)\dbE[X(t)],\dbE[X(t)]\ran \\
\ns\ds\qq\qq\q +\lan R\big\{u(t)-\dbE[u(t)]\big\},u(t)-\dbE[u(t)]\ran+\lan(R+\bar R)\dbE[u(t)],\dbE[u(t)]\ran\Big\}dt \\
\ns\ds\qq\qq \geq \d\dbE\int_0^\infty|u(t)|^2dt,\ea$$
for some $\d>0$. Therefore, under (J) and (S), the map
$u(\cd)\mapsto J(x;u(\cd))$ is a quadratic and coercive functional
on $\cU_{ad}[0,\infty)$. Hence, by a standard argument, we see that
optimal control $u_*(\cd)\in\cU_{ad}[0,\infty)$ must uniquely exist,
and of course, $X_*(\cd)$ is also unique. By a standard argument, we
can show that value function $V(x)$ is of form (\ref{J=Pi}) for some
$\Pi\in\cS^n$, $\Pi\ge0$.

\ms

Now, for any $T>0$, let
$$J_T(x;u(\cd))\1n=\1n\dbE\2n\int_0^T\3n\Big\{\1n\lan QX(t),X(t)\ran\1n+\1n\lan\bar Q\dbE[X(t)],\dbE[X(t)]\ran\1n+\1n\lan\1n Ru(t),u(t)\ran\1n+\1n\lan\1n\bar R\,\dbE[u(t)],\dbE[u(t)]\ran\Big\}dt.$$
We may pose the following problem.

\ms

\bf Problem (LQ)$_T$. \rm For any $x\in\dbR^n$, find a
$u_T(\cd)\in\cU[0,T]$ such that
$$J_T(x;u_T(\cd))=\inf_{u(\cd)\in\cU[0,T]}J_T(x;u(\cd))\equiv V_T(x).$$

By \cite{Yong 2011}, for Problem (LQ)$_T$, under (J), we have a
unique $u_T(\cd)\in\cU[0,T]$ such that
$$V_T(x)=\inf_{u(\cd)\in\cU[0,T]}J_T(x;u(\cd))=J_T(x;u_T(\cd))=\lan\Pi(0;T)x,x\ran,\qq\forall x\in\dbR^n,$$
where
\bel{Riccati-P(t;T)}\left\{\ba{ll}
\ns\ds\dot P(t;T)+P(t;T)A+A^TP(t;T)+C^TP(t;T)C+Q \\
\ns\ds\q-\big[P(t;T)B\2n+\1n C^T\1n P(t;T)D\big] \big[R\1n+\1n D^T\1n P(t;T)D\big]^{-1}\big[B^T\1n P(t;T)\1n+\1n D^T\1n P(t;T)C\big]\1n=0,\q t\in[0,T], \\
\ns\ds P(T;T)=0,\ea\right.\ee
and
$$\left\{\2n\ba{ll}
\ns\ds\dot\Pi(t;T)+\Pi(t;T)(A\1n+\1n\bar A)\1n+\1n(A\1n+\1n\bar A)^T\Pi(t;T)\1n+\1n(C\1n+\1n\bar C)^T\1n P(t;T)(C\1n+\1n\bar C)\1n+\1n Q\1n+\1n\bar Q \\
\ns\ds\q-\big[\Pi(t;T)(B\1n+\1n\bar B)\1n+\1n(C\1n+\1n\bar C)^T\1n P(t;T)(D\1n+\1n\bar D)\big]\big[R+\bar R+(D+\bar D)^TP(t;T)(D+\bar D)\big]^{-1} \\
\ns\ds\q\cd\big[(B\1n+\1n\bar B)^T\Pi(t;T)\1n+\1n(D\1n+\1n\bar D)^T\1n P(t;T)(C\1n+\1n\bar C)\big]=0,\qq t\in[0,T], \\
\ns\ds\Pi(T;T)=0.\ea\right.$$
Further, if we define
$$\left\{\2n\ba{ll}
\ns\ds\G(t;T) = -\big[R+D^TP(t;T)D\big]^{-1}\big[B^TP(t;T)+D^TP(t;T)C\big], \\
\ns\ds\bar\G(t;T) = -\1n\big[R\1n+\1n\bar R\1n+\1n(D\1n+\1n\bar D)^T\1n P(t;T)(D\1n+\1n\bar D)\big]^{-1}\big[(B\1n+\1n\bar B)^T\1n P(t;T)\1n+\1n(D\1n+\1n\bar D)^T\1n P(t;T)(C\1n+\1n\bar C)\big],\ea\right.$$
then the optimal control $u_T(\cd)$ admits the following state
feedback representation:
$$u_T(t)=\G(t;T)\big\{X_T(t)-\dbE[X_T(t)]\big\}+\bar\G(t;T)\dbE[X_T(t)],\qq t\in[0,T],$$
where $X_T(\cd)$ is the solution to the closed-loop system:
\bel{5.38}\left\{\2n\ba{ll}
\ns\ds dX_T(t) = \Big\{\big[A+B\G(t;T)\big]X_T(t)+\big[\bar A+\bar B\bar\G(t;T)+B\big(\bar\G(t;T)-\G(t;T)\big)\big]\dbE[X_T(t)]\Big\}dt \\
\ns\ds\qq\qq\q\2n +\Big\{\big[C+D\G(t;T)\big]X_T(t)+\big[\bar C+\bar D\bar\G(t;T)+D\big(\bar\G(t;T)-\G(t;T)\big)\big]\dbE[X_T(t)]\Big\}dW(t), \\
\ns\ds\qq\qq\qq\qq\qq\qq\qq\qq\qq\qq\qq\qq\qq t\in[0,T], \\
\ns\ds X_T(0) = x.\ea\right.\ee
Observe that (\ref{Riccati-P(t;T)}) coincides with (\ref{Riccati-T0}).
By the proof of Theorem \ref{thm5.1}, we see that
$$\lim_{T\to\infty}P(t;T)=P,\qq t\ge0.$$
Hence,
$$\lim_{T\to\infty}\G(t;T)=-(R+D^TPD)^{-1}(B^TP+D^TPC)\equiv\G,\qq t\ge0.$$
Now, we introduce the following differential Riccati equation (on $[0,\infty)$):
%
$$\left\{\2n\ba{ll}
\ns\ds-\dot{\bar\Pi}(s)+\bar\Pi(s)(A\1n+\1n\bar A)\1n+\1n(A\1n+\1n\bar A)^T\bar\Pi(s)\1n+\1n(C\1n+\1n\bar C)^T\1n\bar P(s)(C\1n+\1n\bar C)\1n+\1n Q\1n+\1n\bar Q\\
\ns\ds\q-\big[\bar\Pi(s)(B\1n+\1n\bar B)\1n+\1n(C\1n+\1n\bar C)^T\1n \bar P(s)(D\1n+\1n\bar D)\big]\big[R+\bar R+(D+\bar D)^T\bar P(s)(D+\bar D)\big]^{-1}\\
\ns\ds\q\cd\big[(B\1n+\1n\bar B)^T\bar\Pi(s)\1n+\1n(D\1n+\1n\bar D)^T\1n\bar P(s)(C\1n+\1n\bar C)\big]=0,\qq t\ge0, \\
\ns\ds\bar\Pi(0)=0.\ea\right.$$
For any $T>0$, let
$$\wt\Pi(t;T)=\bar\Pi(T-t),\qq t\in[0,T].$$
Then by the uniqueness, we must have
$$\Pi(t;T)=\wt\Pi(t;T)=\bar\Pi(T-t),\qq\qq t\in[0,T].$$
Hence,
$$\Pi(0;T)=\bar\Pi(T),\qq\qq T\ge0.$$
Similar to the proof of Theorem \ref{thm5.1}, we have that
$$0\le\Pi(0;T)\le\Pi(0;\bar T)\le\Pi,\qq0\le T\le\bar T<\infty.$$
Thus,
$$\lim_{t\to\infty}\bar\Pi(t)=\bar\Pi\le\Pi.$$
Further, $\bar\Pi$ must be a solution to the following ARE:
$$\ba{ll}
\ns\ds\bar\Pi(A+\bar A)+(A+\bar A)^T\bar\Pi+(C+\bar C)^T\1n P(C+\bar C)+Q+\bar Q \\
\ns\ds\q-\big[\bar\Pi(B+\bar B)+(C+\bar C)^TP(D+\bar D)\big]\big[R+\bar R+(D+\bar D)^TP(D+\bar D)\big]^{-1} \\
\ns\ds\q\cd\big[(B\1n+\1n\bar B)^T\bar\Pi\1n+\1n(D\1n+\1n\bar D)^T\1n P(C\1n+\1n\bar C)\big]=0,\ea$$
Also,
$$\lim_{T\to\infty}\Pi(t;T)=\lim_{T\to\infty}\bar\Pi(T-t)=\bar\Pi,\qq t\ge0.$$
Then
%
$$\lim_{T\to\infty}\bar\G(t;T)=\bar\G_0=-\big[R+\bar R+(D+\bar D)^TP(D+\bar D)\big]^{-1}\big[(B+\bar B)^T\bar\Pi+(D+\bar D)^TP(C+\bar C)\big],\q\forall t\ge0.$$
Recall that $X_T(\cd)$ satisfies (\ref{5.38}). Thus, one has
$$\lim_{T\to\infty}X_T(t;T)=\bar X(t),\qq t\ge0,$$
with $\bar X(\cd)$ being the solution to the following:
$$\left\{\2n\ba{ll}
\ns\ds d\bar X(t) = \ds\Big\{(A+B\G)\bar X(t)+\big[\bar A+\bar B\bar\G_0+B(\bar G_0-\G)\big]\dbE[\bar X(t)]\Big\}dt \\
\ns\ds\qq\qq +\Big\{(C+D\G)\bar X(t)+\big[\bar C+\bar D\bar\G_0+D(\bar \G_0-\G)\big]\dbE[\bar X(t)]\Big\}dW(t),\qq t\in[0,T], \\
\ns\ds\bar X(0) = x.\ea\right.$$
On the other hand,
$$\ba{ll}
\ns\ds \lan\Pi(0;T)x,x\ran=J_T(x;u_T(\cd)) \\
\ns\ds = \dbE\int_0^T\Big\{\lan QX_T(t),X_T(t)\ran+\lan\bar Q\dbE[X_T(t)],\dbE[X_T(t)]\ran \\
\ns\ds\q +\1n\lan\1n R\(\G(t;T)\big\{X_T(t)\1n-\1n\dbE[X_T(t)]\big\}\2n+\1n\bar\G(t;T)\dbE[X_T(t)]\),\G(t;T)\big\{X_T(t)\1n-\1n\dbE[X_T(t)]\}\2n+\1n\bar\G(t;T)\dbE[X_T(t)]\ran \\
\ns\ds\q +\lan\bar R\bar\G(t;T)\dbE[X_T(t)],\bar\G(t;T)\dbE[X_T(t)]\ran\Big\}dt.\ea$$
Thus, sending $T\to\infty$, by Fatou's Lemma, we obtain
$$\ba{ll}
\ns\ds\lan\Pi x,x\ran \geq \lan\bar\Pi x,x\ran\ge\dbE\ds\int_0^\infty\Big\{\lan Q\bar X(t),\bar X(t)\ran+\lan\bar Q\dbE[\bar X(t)],\dbE[\bar X(t)]\ran \\
\ns\ds\qq\qq\q +\1n\lan\1n R\big(\G\big\{\bar X(t)\1n-\1n\dbE[\bar X(t)]\big\}\2n+\1n\bar\G_0\dbE[\bar X(t)]\big),\G\big\{\bar X(t)\1n-\1n\dbE[\bar X(t)]\}\2n+\1n\bar\G_0\dbE[\bar X(t)]\ran\\
\ns\ds\qq\qq\q +\lan\bar R\bar\G_0\dbE[\bar X(t)],\bar\G_0\dbE[\bar X(t)]\ran\Big\}dt=J(x;\bar u(\cd))\ge\lan\Pi x,x\ran.\ea$$
Hence,
$$\bar\Pi=\Pi,\qq\bar\G_0=\bar\G,$$
and $(\G,\bar\G)$ is an MF-$L^2_{Q,\bar Q}$-stabilizer of the
system, and $(\bar X(\cd),\bar u(\cd))=(X_*(\cd),u_*(\cd))$ is the
optimal pair. \endpf

\section{Optimal MF-LQ Controls Presented via Tackling AREs}

\subsection{Tackling AREs via LMIs}

One of the main ideas of this section is to reformulate the AREs as
{\it linear matrix inequalities} (LMIs, for short). Let us introduce
the general notion of LMIs according to \cite{AitZhou,LZAit}, and develop it to solve our mean-field LQ problem.

\bde \rm
Let $ F_0,F_1,\cdots,F_m \in {\cal S}^n$ be given.
Inequalities consisting of any combination of  the following relations
\bel{LMI-def}\ba{ll}
\ns\ds F(x) \deq F_0 + \ds\sum_{i=1}^m x_iF_i > 0, \qq\mbox{ or }\qq F(x) \deq F_0 + \ds\sum_{i=1}^m x_iF_i \geq 0,\ea\ee
are called LMIs with respect to the variable $x=(x_1,\cdots,x_m)^T \in \dbR^m$.
When the LMI is satisfied by a vector $x$ we say  that
the LMI is feasible and $x$ is a feasible point.
\ede

Next, let us state some facts about general {\it semi-definite programming} (SDP, for short) problems and their duals.

\bde \rm
Let $c=(c_1,\cdots,c_m)^T\in\dbR^m$ and $F_0,F_1,\ldots,F_m \in {\cal S}^n$ be given. The following optimization problem
\bel{SDP-def}\ba{rl}
\min & c^Tx, \\
\mbox{\rm subject to} & F(x)\equiv F_0 + \ds\sum_{i=1}^m x_iF_i  \geq 0, \ea\ee
is called a semidefinite programming. Moreover,
the dual problem of  the SDP (\ref{SDP-def}) is defined as
\bel{Dual-def}\ba{rl}
\max & -{\bf Tr}(F_0Z), \\
\mbox{\rm subject to} & Z\in {\cal S}^n,\;\;  {\bf Tr}(ZF_i)=c_i, \;\;i=1,2,\cdots,m, \;\; Z\geq 0. \ea\ee
\ede

\ms

The following basic assumption is imposed throughout this section.

\begin{assumption} \rm
The controlled MF-FSDE system $[A,\bar A,C,\bar C;B,\bar B,D,\bar D]$ is MF-$L^2$-stabilizable.
\end{assumption}

\ms

For notational convenience, we rewrite the AREs (\ref{A-Riccati}) as follows
\bel{ric-QR}
\ds \cR(P,Q,\Q,R,\R) = 0, \qq \bar\cR(P,\Pi,Q,\Q,R,\R) = 0, \ee
where
$$\left\{\2n\ba{ll}
\ns\ds \cR(P,Q,\Q,R,\R)\deq PA+A^TP+C^TPC-(PB+C^TPD)(R+D^TPD)^{-1}(B^TP+D^TPC)+Q, \\
\ns\ds \bar\cR(P,\Pi,Q,\Q,R,\R)\deq\Pi(A + \A) + (A + \A)^T\Pi + (C + \C)^TP(C + \C) + Q + \Q \\
\ns\ds\qq\qq\qq\qq\qq -\big[\Pi(B + \B) + (C + \C)^TP(D + \D)\big]\big[R + \R + (D + \D)^TP(D + \D)\big]^{-1} \\
\ns\ds\qq\qq\qq\qq\qq \cdot\big[(B + \B)^T\Pi + (D + \D)^TP(C + \C)\big]. \ea\right.$$

\ms

\bl\label{lem-increase} \sl
Let $Q_1,\Q_1,Q_2,\Q_2\in\cS^n$ and $R_1,\R_1,R_2,\R_2\in\cS^m$ be given satisfying
$$Q_1\le Q_2,\q\Q_1\le\Q_2,\q R_1\le R_2,\q\R_1\le\R_2.$$
Assume that there exists $(P_0,\Pi_0)$ such that
$$\cR(P_0,Q_1,\Q_1,R_1,\R_1)>0, \qq \bar\cR(P_0,\Pi_0,Q_1,\Q_1,R_1,\R_1)>0.$$
Then there exist $(P_1^*,\Pi_1^*)$ and $(P_2^*,\Pi_2^*)$ satisfying
$$\left\{\2n\ba{ll}
\ds \cR(P_i^*,Q_i,\Q_i,R_i,\R_i)=0, \qq \bar\cR(P_i^*,\Pi_i^*,Q_i,\Q_i,R_i,\R_i)=0, \qq\mbox{for } i = 1, 2, \\
\ns\ds P_1^* \leq P_2^* \quad\mbox{ and }\quad \Pi_1^* \leq \Pi_2^*. \ea\right.$$
\el

\it Proof. \rm By the assumptions of this Lemma, $(P_0,\Pi_0)$ must
also satisfy
$$\cR(P_0,Q_2,\Q_2,R_2,\R_2)>0, \qq \bar\cR(P_0,\Pi_0,Q_2,\Q_2,R_2,\R_2)>0.$$
It then follows from Proposition \ref{coro-sol1} that there exist $(P_1^*,\Pi_1^*)$ and $(P_2^*,\Pi_2^*)$, which are the maximal solutions of their respective AREs:
$$\cR(P_i^*,Q_i,\Q_i,R_i,\R_i)=0, \qq \bar\cR(P_i^*,\Pi_i^*,Q_i,\Q_i,R_i,\R_i)=0, \qq\mbox{for } i = 1, 2.$$
Furthermore, $(P_1^*,\Pi_1^*)$ must satisfy
$$\cR(P_1^*,Q_2,\Q_2,R_2,\R_2)\ge0,\qq\bar\cR(P_1^*,\Pi_1^*,Q_2,\Q_2,R_2,\R_2)\ge0.$$
Hence $P_1^* \leq P_2^*$ and $\Pi_1^* \leq \Pi_2^*$ because
$(P_2^*,\Pi_2^*)$ is the maximal solution to its AREs.
\endpf

\ms

Consider the following SDP problem
\bel{ultime-SDP-Pb}\ba{rl}
\ns\ds \max & {\bf Tr}(P) + {\bf Tr}(\Pi), \\
\ns\ds \mbox{\rm subject to} & \left\{\2n\ba{ll}\left[\ba{c|c} PA+A^TP+C^TPC+Q & PB+C^TPD \\ \hline B^TP+D^TPC & R+D^TPD \ea\right] \geq 0, \\
\ns\ds \left[\ba{c|c}\ba{ll}\Pi(A + \A) + (A + \A)^T\Pi \\ + (C + \C)^TP(C + \C) + Q + \Q \ea & \Pi(B + \B)+(C + \C)^TP(D + \D) \\ \hline
(B + \B)^T\Pi+(D + \D)^TP(C + \C)  & R+\R+(D+\D)^TP(D+\D) \ea\right]\geq 0.
%
%
\ea\right.\ea\ee

\ms

\bt\label{ultime-theo} \sl
Let $Q,\Q\in\cS^n$, $R,\R\in\cS^m$ be given. The following are equivalent:
\begin{itemize}
\item[{\rm (i)}]  There exists $(P_0,\Pi_0)$ such that $\cR(P_0,Q,\Q,R,\R)\geq 0$ and $\bar\cR(P_0,\Pi_0,Q,\Q,R,\R)\geq 0$.
\item[{\rm (ii)}] There  exists a solution to the AREs $(\ref{ric-QR})$.
\end{itemize}
Moreover, when (i) or (ii) holds, the AREs $(\ref{ric-QR})$ has a
maximal solution $(P^*,\Pi^*)$ which is the unique optimal solution to the SDP problem $(\ref{ultime-SDP-Pb})$.
\et

\it Proof. \rm We only need to prove that (i) implies (ii). Let
$P_0$ be given as in (i). For any $ \epsilon >0$ and $\bar\epsilon >
0$, we have $\cR(P_0,Q+\epsilon I,Q+\bar\epsilon I,R,\R)> 0$ and
$\bar\cR(P_0,\Pi_0,Q+\epsilon I,Q+\bar\epsilon I,R,\R)> 0$. Applying
Proposition \ref{coro-sol1} and Lemma \ref{lem-increase}, we have
that for  any  positive decreasing  sequence $\epsilon_i \rightarrow
0$ and $\bar\epsilon_i \rightarrow 0$ there exists a decreasing
sequence of symmetric matrices
$$P_{\epsilon_0} \geq \cdots \geq P_{\epsilon_i}\geq  P_{\epsilon_{i+1}} \geq P_0, \qq \Pi_{\bar\epsilon_0} \geq \cdots \geq \Pi_{\bar\epsilon_i} \geq  \Pi_{\bar\epsilon_{i+1}} \geq \Pi_0$$
such that $\cR(P_{\epsilon_i},Q+\epsilon_i I,\Q+\bar\epsilon_i
I,R,\R)=0$ and $\bar\cR(P_{\epsilon_i},\Pi_{\bar\epsilon_i},\Q+\epsilon_i I,Q+\bar\epsilon_i I,R,\R)=0$.
Hence the limit $P^*=\ds\lim_{\epsilon_i\rightarrow 0}P_{\epsilon_i}$ and
$\Pi^*=\ds\lim_{\bar\epsilon_i\rightarrow 0}\Pi_{\bar\epsilon_i}$ exist and satisfy
$$\cR(P^*,Q,\Q,R,\R) = 0, \qq \bar\cR(P^*,\Pi^*,Q,\Q,R,\R) = 0.$$
In addition, $(P^*,\Pi^*)$ must be the maximal solution of  the AREs
due to the arbitrariness of $(P_0,\Pi_0)$. By Schur's lemma (Lemma \ref{lem-schur}), $(P^*,\Pi^*)$
is an optimal solution to the problem (\ref{ultime-SDP-Pb}) due to its maximality.
To prove the uniqueness, let $(P_*,\Pi_*)$ be any  optimal solution to (\ref{ultime-SDP-Pb}).
Then ${\bf Tr}(P^*-P_*)+{\bf Tr}(\Pi^*-\Pi_*)=0$ as both $(P^*,\Pi^*)$ and $(P_*,\Pi_*)$ are optimal to (\ref{ultime-SDP-Pb}). However,
$P^*-P_*\geq 0$ and $\Pi^*-\Pi_*\geq 0$ since $(P^*,\Pi^*)$ is the maximal solution of (\ref{ultime-SDP-Pb}).
This yields $P^*-P_*=0$ and $\Pi^*-\Pi_*=0$.
\endpf

\ms

As an  immediate consequence of Theorem \ref{ultime-theo}, we have the following result for the standard case $Q, \Q \geq  0$ and $R, \R > 0$.

\bc \sl
If $Q, \Q \geq  0$ and $R, \R > 0$, then the AREs $(\ref{ric-QR})$ admits a maximal solution $(P^*,\Pi^*)$ with
$P^*,\Pi^* \geq 0$ which is also the unique solution to the SDP
$(\ref{ultime-SDP-Pb})$. In addition, if $Q,\Q >  0$ and $R,\R> 0$,
then the maximal solution $(P^*,\Pi^*)$ with $P^*,\Pi^* > 0$ and the feedback control
$$\ba{ll}
\ns\ds\3n   u^*(t) = -(R+D^TP^*D)^{-1}(B^TP^*+D^TP^*C)\big(X^*(t)-\dbE[X^*(t)]\big) \\
\ns\ds\qq\q -\big(R+\R+(D+\D)^TP^*(D+\D)\big)^{-1}\big[(B + \B)^T\Pi^*+(D + \D)^TP^*(C + \C)\big]\dbE[X^*(t)]\ea$$
is stabilizing for the system (\ref{MF-state1}).
\ec

\it Proof.
\rm
When $Q,\Q \geq  0$ and $R,\R > 0$, $(P_0,\Pi_0)=(0,0)$ satisfies the LMIs
\bel{aux}
\ds \cR(P,Q,\Q,R,\R)\geq 0, \qq \bar\cR(P,\Pi,Q,\Q,R,\R)\geq 0. \ee
Hence by Theorems \ref{ultime-theo} the  AREs (\ref{ric-QR})
admits a maximal solution  $(P^*,\Pi^*)$. Moreover, by the proof of
Theorems \ref{ultime-theo}, $P^* \geq P_0=0$ and $\Pi^* \geq \Pi_0=0$.
If in addition $Q,\Q >  0$ and $R,\R > 0$, then $(\tilde P_0,\tilde\Pi_0) = (\delta I,\bar\delta I)$
solves (\ref{aux}) for a sufficiently small $\delta,\bar\delta > 0$.
Hence $P^*\geq \tilde P_0=\delta I>0$ and $\Pi^*\geq \tilde\Pi_0=\bar\delta I>0$.
Moreover, by virtue of Proposition \ref{theo-stab-sol}, the corresponding feedback control is stabilizing since (\ref{aux}) is
strictly feasible in this case.
\endpf

\subsection{Optimal feedback Control}

In this subsection, we show that
the value function of Problem MF-LQ
can be expressed in terms of the maximal
solution to the AREs (\ref{ric-QR}). Moreover,
if there exists an  optimal control  of Problem MF-LQ  then it is
necessarily represented as a feedback via the  maximal solution to the AREs.

\bt\label{theo-vcost=Pmax} \sl
Assume that Theorem \ref{ultime-theo}-(i) holds. Then Problem (MF-LQ) is well-posed and
the value function is given by  $V(x)=x^T\Pi^*x$, $\forall x \in
\dbR^n$, where $(P^*,\Pi^*)$ is the maximal solution to the AREs
$(\ref{ric-QR})$. \et

\it Proof.
\rm
The well-posedness has been shown in Theorem \ref{thm5.2}, which also yields $V(x) = x^T\Pi^*x$.

Now, for any fixed $\epsilon >0$, the LMIs
\bel{auxlmi}
\ds \cR(P,Q+\epsilon I,\Q+\epsilon I,R,\R) \geq 0, \qq \bar\cR(P,\Pi,Q+\epsilon I,\Q+\epsilon I,R,\R) \geq 0 \ee
are strictly feasible. Hence by Proposition \ref{coro-sol1}, there is
a maximal solution, denoted by $(P_\epsilon,\Pi_{\epsilon})$, to the corresponding AREs
$$\cR(P,Q+\epsilon I,\Q+\epsilon I,R,\R) = 0, \qq \bar\cR(P,\Pi,Q+\epsilon I,\Q+\epsilon I,R,\R) = 0.$$
In addition, by Proposition \ref{theo-stab-sol}, the feedback control
$u_{\epsilon}(t)=\Gamma_{\epsilon}\big(X_{\epsilon}(t)-\dbE[X_\epsilon(t)]\big)+\bar\Gamma_\epsilon\dbE[X_\epsilon(t)]$
is stabilizing, where
$$\left\{\2n\ba{ll}
\ns\ds \Gamma_{\epsilon}=-(R+D^TP_{\epsilon}D)^{-1}(B^TP_{\epsilon}+D^TP_{\epsilon}C), \\
\ns\ds \bar\Gamma_{\epsilon}=-\big(R+\R+(D+\D)^TP_{\epsilon}(D+\D)\big)^{-1}\big[(B + \B)^T\Pi_{\epsilon}+(D + \D)^TP_{\epsilon}(C + \C)\big].\ea\right.$$
It is easy to verify that $P_{\epsilon},\Pi_{\epsilon},\Gamma_{\epsilon}$ and $\bar\Gamma_{\epsilon}$ satisfy the following equations
\bel{lyap-epsilon}\left\{\2n\ba{ll}
\ns\ds (A+B\Gamma_{\epsilon})^TP_{\epsilon}+P_{\epsilon}(A+B\Gamma_{\epsilon})+(C+DK_{\epsilon})^TP_{\epsilon}(C+D\Gamma_{\epsilon}) = -Q-\epsilon I-\Gamma_{\epsilon}^TR\Gamma_{\epsilon}, \\
\ns\ds (A+\A+B\bar\Gamma_{\epsilon}+\B\bar\Gamma_{\epsilon})\Pi_{\epsilon} + \Pi_{\epsilon}(A+\A+B\bar\Gamma_{\epsilon}+\B\bar\Gamma_{\epsilon})^T \\
\ns\ds\q + (C+\C+D\bar\Gamma_{\epsilon}+ \D\bar\Gamma_{\epsilon})^TP_{\epsilon}(C+\C+D\bar\Gamma_{\epsilon}+ \D\bar\Gamma_{\epsilon}) = -Q-\Q-2\epsilon I-\bar\Gamma_{\epsilon}^T(R+\R)\bar\Gamma_{\epsilon}. \ea\right.\ee
Applying Lemma \ref{lem-equality} to $M=P_{\epsilon}$, $N=\Pi_{\epsilon}$ and substituting $u_{\epsilon}(t)$ into (\ref{equa-M}), we have
$$\ba{ll}
\ns\ds \dbE\int_0^t\Big\{\langle(Q+\epsilon I)X_{\epsilon}(s),X_{\epsilon}(s)\rangle + \langle(\Q+\epsilon I)\dbE[X_{\epsilon}(s)],\dbE[X_{\epsilon}(s)]\rangle\ds \\
\ns\ds\qq + \langle Ru_{\epsilon}(s),u_{\epsilon}(s)\rangle + \langle\R\dbE[u_{\epsilon}(s)],\dbE[u_{\epsilon}(s)]\rangle\Big\}ds \\
\ns\ds = \dbE\int_0^t\Big\{\langle(Q+\epsilon I)\big(X_{\epsilon}(s)-\dbE[X_{\epsilon}(s)]\big),\big(X_{\epsilon}(s)-\dbE[X_{\epsilon}(s)]\big)\rangle + \langle(Q+\Q+2\epsilon I)\dbE[X_{\epsilon}(s)],\dbE[X_{\epsilon}(s)]\rangle \\
\ns\ds\qq + \langle R\big(u_{\epsilon}(s)-\dbE[u_{\epsilon}(s)]\big),\big(u_{\epsilon}(s)-\dbE[u_{\epsilon}(s)]\big)\rangle + \langle(R+\R)\dbE[u_{\epsilon}(s)],\dbE[u_{\epsilon}(s)]\rangle\Big\}ds \\
\ns\ds = \dbE\int_0^t\Big\{\langle (Q+\epsilon I+\Gamma_\epsilon^TR\Gamma_\epsilon)\big(X_{\epsilon}(s)-\dbE[X_{\epsilon}(s)]\big),\big(X_{\epsilon}(s)-\dbE[X_{\epsilon}(s)]\big)\rangle \\
\ns\ds\qq + \langle\big(Q+\Q+2\epsilon I+\bar\Gamma_\epsilon^T(R+\R)\bar\Gamma_\epsilon\big)\dbE[X_{\epsilon}(s)],\dbE[X_{\epsilon}(s)]\rangle\Big\}ds \\
\ns\ds \equiv \dbE\int_0^t\Big\{\big(X_{\epsilon}(s)-\dbE[X_{\epsilon}(s)]\big)^T\big[(A+B\Gamma_\epsilon)^TP_\epsilon+P_\epsilon(A+B\Gamma_\epsilon)+(C+D\Gamma_\epsilon)^TP_\epsilon(C+D\Gamma_\epsilon)\big]\big(X_{\epsilon}(s)-\dbE[X_{\epsilon}(s)]\big) \\
\ns\ds\qq + \dbE[X_{\epsilon}(s)]^T\big[(C+\C+D\bar\Gamma_\epsilon+\D\bar\Gamma_\epsilon)^TP_\epsilon(C+\C+D\bar\Gamma_\epsilon+\D\bar\Gamma_\epsilon) \\
\ns\ds\qq + (A+\A+B\bar\Gamma_\epsilon+\B\bar\Gamma_\epsilon)\Pi_\epsilon + \Pi_\epsilon(A+\A+B\bar\Gamma_\epsilon+\B\bar\Gamma_\epsilon)^T\big]\dbE[X_{\epsilon}(s)]\Big\}ds \\
\ns\ds = -\dbE\Big[\big(X_{\epsilon}(t)-\dbE[X_{\epsilon}(t)]\big)^TP_\epsilon\big(X_{\epsilon}(t)-\dbE[X_{\epsilon}(t)]\big)\Big] + x^T\Pi_\epsilon x - \dbE[X_{\epsilon}(t)]^T\Pi_\epsilon\dbE[X_{\epsilon}(t)]. \ea$$
Since
$\ds\lim_{t\rightarrow +\infty}\dbE\big[\big(X_{\epsilon}(t)-\dbE[X_{\epsilon}(t)]\big)^TP_\epsilon\big(X_{\epsilon}(t)-\dbE[X_{\epsilon}(t)]\big)\big]=0$
and $\ds\lim_{t\rightarrow +\infty}\dbE[X_{\epsilon}(t)]^T\Pi_\epsilon\dbE[X_{\epsilon}(t)]=0$, we obtain
$$\ba{ll}
\ns\ds x^T\Pi_\epsilon x = \dbE\int_0^{\infty}\Big\{\langle(Q+\epsilon I)X_{\epsilon}(s),X_{\epsilon}(s)\rangle + \langle(\Q+\epsilon I)\dbE[X_{\epsilon}(s)],\dbE[X_{\epsilon}(s)]\rangle\ds \\
\ns\ds\qq\qq\q + \langle Ru_{\epsilon}(s),u_{\epsilon}(s)\rangle + \langle\R\dbE[u_{\epsilon}(s)],\dbE[u_{\epsilon}(s)]\rangle\Big\}ds \geq V(x). \ea$$
On the other hand, since $P^*=\ds\lim_{\epsilon \rightarrow 0}P_\epsilon$ and $\Pi^*=\ds\lim_{\epsilon \rightarrow 0}\Pi_\epsilon$
(similar to the proof of Theorem \ref{ultime-theo}), we have $V(x) \leq  x^T\Pi^*x $.
This completes the proof.
\endpf

\ms

\bc
\sl
Assume that Theorem \ref{ultime-theo}-(i) holds. If there exists an optimal control
of Problem (MF-LQ), then it must be unique and represented
by the state feedback control
$$u^*(t)=\Gamma^*\big(X^*(t)-\dbE[X^*(t)]\big)+\bar\Gamma^*\dbE[X^*(t)],$$
where $(P^*,\Pi^*)$ is the maximal solution to the AREs
$(\ref{ric-QR})$, and
$$\left\{\2n\ba{ll}
\ns\ds \Gamma^* = -(R+D^TP^*D)^{-1}(B^TP^*+D^TP^*C), \\
\ns\ds \bar\Gamma^* = -\big(R+\R+(D+\D)^TP^*(D+\D)\big)^{-1}\big[(B + \B)^T\Pi^*+(D + \D)^TP^*(C + \C)\big]. \ea\right.$$
\ec

\it Proof.
\rm
Let $(X^*(\cdot),u^*(\cdot))$  be an optimal pair of
the LQ problem. Then a completion of squares shows
$$\ba{ll}
\ns\ds \dbE\int_0^t\Big\{\langle QX^*(s),X^*(s)\rangle+\langle\Q\dbE[X^*(s)],\dbE[X^*(s)]\rangle+\langle Ru^*(s),u^*(s)\rangle+\langle\R\dbE[u^*(s)],\dbE[u^*(s)]\rangle\Big\}ds \\
\ns\ds = \dbE\int_0^t\Big\{\langle Q\big(X^*(s)-\dbE[X^*(s)]\big),\big(X^*(s)-\dbE[X^*(s)]\big)\rangle + \langle(Q+\Q)\dbE[X^*(s)],\dbE[X^*(s)]\rangle \\
\ns\ds\qq + \langle R\big(u^*(s)-\dbE[u^*(s)]\big),\big(u^*(s)-\dbE[u^*(s)]\big)\rangle + \langle(R+\R)\dbE[u^*(s)],\dbE[u^*(s)]\rangle\Big\}ds \\
\ns\ds\qq -\dbE\Big[\big(X^*(t)-\dbE[X^*(t)]\big)^TP^*\big(X^*(t)-\dbE[X^*(t)]\big)\Big] + x^T\Pi^*x-\dbE[X^*(t)]^T\Pi^*\dbE[X^*(t)] \\
\ns\ds\qq + \dbE\int_0^t\Big\{u^*(s)-\dbE[u^*(s)]-\Gamma^*\big(X^*(s)-\dbE[X^*(s)]\big)\Big]^T(R+D^TP^*D)^{-1} \\
\ns\ds\qq \cdot\Big[u^*(s)-\dbE[u^*(s)]-\Gamma^*\big(X^*(s)-\dbE[X^*(s)]\big)\Big\}ds \\
\ns\ds\qq + \dbE\int_0^t\Big\{\dbE[u^*(s)]-\bar\Gamma^*\dbE[X^*(s)]\Big]^T\big(R+\R+(D+\D)^TP^*(D+\D)\big)^{-1}\Big[\dbE[u^*(s)]-\bar\Gamma^*\dbE[X^*(s)]\Big\}ds. \ea$$
As $u^*(\cdot)$ is stabilizing, we have
$$\lim_{t\rightarrow +\infty}\dbE\big[\big(X^*(t)-\dbE[X^*(t)]\big)^TP^*\big(X^*(t)-\dbE[X^*(t)]\big)\big]=0, \q
\lim_{t\rightarrow +\infty}\dbE[X^*(t)]^T\Pi^*\dbE[X^*(t)]=0,$$
which implies
\bel{optimality}\ba{ll}
\ns\ds V(x) = J(x,u^*(\cdot)) \\
\ns\ds = x^T\Pi^*x + \ds\dbE\int_0^{\infty}\Big\{u^*(s)-\dbE[u^*(s)]-\Gamma^*\big(X^*(s)-\dbE[X^*(s)]\big)\Big]^T(R+D^TP^*D)^{-1} \\
\ns\ds\q \cdot\Big[u^*(s)-\dbE[u^*(s)]-\Gamma^*\big(X^*(s)-\dbE[X^*(s)]\big)\Big\}ds \\
\ns\ds\q + \dbE\int_0^{\infty}\Big\{\dbE[u^*(s)]-\bar\Gamma^*\dbE[X^*(s)]\Big]^T\big(R+\R+(D+\D)^TP^*(D+\D)\big)^{-1}\Big[\dbE[u^*(s)]-\bar\Gamma^*\dbE[X^*(s)]\Big\}ds. \ea\ee
By Theorem \ref{theo-vcost=Pmax} we have  $V(x)=x^T\Pi^*x$. Hence,
$$\left\{\2n\ba{ll}
\ns\ds \dbE\int_0^{\infty}\Big\{u^*(s)-\dbE[u^*(s)]-\Gamma^*\big(X^*(s)-\dbE[X^*(s)]\big)\Big]^T(R+D^TP^*D)^{-1} \\
\ns\ds\q \cdot\Big[u^*(s)-\dbE[u^*(s)]-\Gamma^*\big(X^*(s)-\dbE[X^*(s)]\big)\Big\}ds = 0, \\
\ns\ds \dbE\int_0^{\infty}\Big\{\dbE[u^*(s)]-\bar\Gamma^*\dbE[X^*(s)]\Big]^T\big(R+\R+(D+\D)^TP^*(D+\D)\big)^{-1}\Big[\dbE[u^*(s)]-\bar\Gamma^*\dbE[X^*(s)]\Big\}ds = 0. \ea\right.$$
As $R+D^TP^*D$ and $R+\R+(D+\D)^TP^*(D+\D)$ are constant positive definite matrices,
$u^*(t)$ has to be in a feedback form
$u^*(t)=\Gamma^*\big(X^*(t)-\dbE[X^*(t)]\big)+\bar\Gamma^*\dbE[X^*(t)]$.
\endpf

\section{Numerical Examples}
In this section, we report our numerical
experiments based on the approach developed in the previous sections.
Note that the numerical  algorithm we have used
for checking LMIs or solving SDP \cite{VaB:96}.

The system dynamics (\ref{MF-state1}) in our experiments is specified by the following matrices
$$
\ba{ll}
A = \left[\ba{rrrrr}
    -0.7 & 0.2 & -0.9 & -0.2 & -0.7 \\
    1.0 & -0.6 & 2.0 & -0.6 & -0.8 \\
    0.8 & 0.8 & -1.7 & -1.5 & 1.1 \\
    0.7 & -0.2 & 0.1 & -0.3 & -0.2 \\
    0.6 & -1.0 & -1.3 & 0.6 & -0.2
\ea\right],
&
\bar A = \left[\ba{rrrrr}
    -0.75 & 0.25 & -0.95 & -0.25 & -0.75 \\
    1.05 & -0.65 & 2.05 & -0.65 & -0.85 \\
    0.85 & 0.85  & -1.75 & -1.55 & 1.15 \\
    0.75 & -0.25 & 0.15 & -0.35 & -0.25 \\
    0.65 & -1.05 & -1.35 & 0.65 & -0.25
\ea\right], \\ [12mm]
B = \left[\ba{rr}
    1.4 & -0.7 \\
    0.3 & -1.7 \\
    0.1 & -1.7 \\
    -0.1 & 0.1 \\
    0.4 & -1.2
\ea\right],
&
\bar B = \left[\ba{rr}
    1.45 & -0.75 \\
    0.35 & -1.75 \\
    0.15 & -1.75 \\
    -0.15 & 0.15 \\
    0.45 & -1.25 \\
\ea\right], \\ [12mm]
C = \left[\ba{rrrrr}
    0.1  & 0.1  & 0.2  & -0.1 &  0.4 \\
    -0.1 & -0.3 & 0.2  & -0.1 & -0.3 \\
    0.6  & 0.4  & -0.3 &  0.1 & -0.2 \\
    -0.1 & 0.2  & -0.2 & -0.1 &  0.1 \\
    -0.2 & 0.2  & 0.3  & 0.2  & -0.3
\ea\right],
&
\bar C = \left[\ba{rrrrr}
    0.15  &  0.15 &  0.25 & -0.15 &  0.45 \\
    -0.15 & -0.35 &  0.25 & -0.15 & -0.35 \\
    0.65  & 0.45 & -0.35  & 0.15  & -0.25 \\
    -0.15 & 0.25 & -0.25  & -0.15 &  0.15 \\
    -0.25 & 0.25 &  0.35  & 0.25  & -0.35
\ea\right], \\ [12mm]
D = \left[\ba{rr}
     0.7 &  -0.3 \\
     0.2 &  -0.8 \\
     0.1 & -0.8 \\
     -0.1 &  0.5 \\
      0.2 & -0.6
\ea\right],
&
\bar D = \left[\ba{rr}
     0.75 & -0.35 \\
     0.25 & -0.85 \\
     0.15 & -0.85 \\
     -0.15 &  0.55 \\
      0.25 & -0.65
\ea\right].
\ea
$$

\subsection{Numerical test of MF-L$^2$ stabilizability}

Since we have shown that the controlled MF-FSDE system is MF-L$^2$-stabilizable
in Proposition \ref{theo-stab} if and only if (\ref{cond-LMI}) is feasible (with respect to
the variables $\mathbb{X}$, $\bar{\mathbb{X}}$, $Y$ and $\bar Y$), we should check the MF-L$^2$ stabilizability first by tackling inequalities.
After running the calculation of SDP program via Matlab software, the obtained feasible matrices $\mathbb{X}$, $\bar{\mathbb{X}}$, $Y$ and $\bar Y$ satisfy Proposition \ref{theo-stab}:
$$
\mathbb{X} = \left[\ba{rrrrr}
    26.1032 &   0.6379  & -7.9410  &  1.4143  & -7.4032 \\
    0.6379  & 17.0911   & -0.4114  &  8.2578  &  1.3415 \\
   -7.9410  & -0.4114   & 19.4946  &  1.1492  & 14.0620 \\
    1.4143  &  8.2578   & 1.1492   & 21.8509  &  7.8151 \\
   -7.4032  &  1.3415   & 14.0620  &  7.8151  & 40.5193
\ea\right],
$$

$$
\bar{\mathbb{X}} = \left[\ba{rrrrr}
    0.0471 &  -0.0617  &  0.0114  & -0.2361  & -0.0333 \\
   -0.0617 &  -0.1398  & -0.1104  &  0.2431  &  0.3623 \\
    0.0114 &  -0.1104  &  0.0283  &  0.1159  &  0.0443 \\
   -0.2361 &   0.2431  &  0.1159  &  0.4583  &  0.0880 \\
   -0.0333 &   0.3623  &  0.0443  &  0.0880  &  0.0952
\ea\right],
$$

$$
Y = \left[\ba{rrrrr}
    -12.1167  & -1.8513  &  7.0876 & -11.3987 &  -1.6418 \\
    0.9756    & 2.1581   & 5.2614  & -16.0940 & -12.8827
\ea\right]
$$
and
$$
\bar Y = \left[\ba{rrrrr}
   -0.3539 &  -0.0281 &  -0.0278 &  -0.2997  &  0.1924 \\
   -0.0070 &  -0.0900 &   0.1334 &  -0.4658  &  0.1065
\ea\right]
$$
which give rise to the stabilizing feedback control law $u(t)=K(X(t) - \dbE[X(t)]) + \bar K\dbE[X(t)]$
with the following feedback gain
$$
K = Y\mathbb{X}^{-1}\left[\ba{rrrrr}
   -0.3725  &  0.1843  &  0.3405  & -0.5390 &  -0.1289 \\
    0.1689  &  0.5864  &  0.6700  & -0.8716 &  -0.3709
\ea\right]
$$
and
$$
\bar K = \bar Y\bar{\mathbb{X}}^{-1} = \left[\ba{rrrrr}
    4.0644  &  1.2449  & -3.7655  &  1.9996  & -1.3936 \\
    4.2782  &  0.3405  &  0.6593  &  0.7858  &  0.2849
\ea\right].
$$

\subsection{Numerical solutions of SARE}

Now we tackle the SARE (\ref{ric-QR})
for the following $Q$, $\bar Q$, $R$ and $\bar R$ via
solving the SDP problem (\ref{ultime-SDP-Pb}):
$$
Q = \mbox{\rm diag}([0,1,1,0,1]) \qq \mbox{ and } \qq \bar Q = \mbox{\rm diag}([0,0.5,1,0,0.5]),
$$
and
$$
R =\mbox{\rm diag}([1,1]) \qq \mbox{ and } \qq \bar R = \mbox{\rm diag}([1.5,1]).
$$
We then gain the following solution $(P,\Pi)$
$$
P = \left[\ba{rrrrr}
    0.4151  &  0.3890  &  0.2068  &  0.0162  & -0.4059 \\
    0.3890  &  2.7208  &  1.9097  & -2.6074  & -0.7756 \\
    0.2068  &  1.9097  &  1.8535  & -1.8330  & -0.8979 \\
    0.0162  & -2.6074  & -1.8330  &  4.2403  & -0.2665 \\
   -0.4059  & -0.7756  & -0.8979  & -0.2665  &  2.1537
\ea\right]
$$
and
$$
\Pi = \left[\ba{rrrrr}
    0.6147  &  0.5721  &  0.2644  & -0.1455  & -0.6138 \\
    0.5721  &  4.2579  &  2.8706  & -4.4158  & -0.6536 \\
    0.2644  &  2.8706  &  2.6758  & -2.6653  & -1.0890 \\
   -0.1455  & -4.4158  & -2.6653  &  6.8158  & -1.0674 \\
   -0.6138  & -0.6536  & -1.0890  & -1.0674  &  3.1641
\ea\right].
$$

\rm

\appendix

\section{Appendix}

\subsection{Some useful lemmas}

The well-known Schur lemma in \cite{BEFB:94} plays a key technical
role.

\begin{lemma}[Schur's lemma]\label{lem-schur} \sl
Let matrices $M=M^T, N$ and  $R=R^T>0$  be given with  appropriate
 dimensions.  Then the following conditions
are equivalent:

{\rm (i)} $M-NR^{-1}N^T \geq \mbox{ (resp. $>$) }0.$

{\rm (ii)} $\left[\ba{cc}
M & N \\
N^T & R
\ea\right] \geq \mbox{ (resp. $>$) }0.$

{\rm (iii)} $\left[\ba{cc}
R & N^T \\
N & M
\ea\right] \geq \mbox{ (resp. $>$) }0.$
\end{lemma}

In the original Schur lemma, the matrix $R$ is required to be
nonsingular. When $R$ is possibly singular, we have an extended
Schur's lemma making use of some generalized inverse matrices. To
make it more precise, for any matrix $M$, there exists a  unique
matrix $M^+ $, called the {\it Moore-Penrose inverse} \cite{Pen:55},
such that
$$MM^+ M=M, \; M^+ MM^+ =M^+, \; (MM^+ )^T=MM^+, \; (M^+ M)^T=M^+ M.$$

\begin{lemma} \label{pseudo-lem} \sl
For a symmetric matrix $S$, we have
\begin{itemize}
\item[{\rm (i)}]  $S^+ = ({S^+})^T$.
\item[{\rm (ii)}] $S\geq 0$  if and only if $ S^+ \geq 0$.
\item[{\rm (iii)}] $SS^+ =S^+ S$.
\end{itemize}
\end{lemma}

Its proof can be found in \cite{AitZhou}.

\begin{lemma}[Extended Schur's lemma]\label{schur-ext} \sl
Let matrices $M=M^T, N$ and  $R=R^T$ be given with
appropriate dimensions. Then the following conditions are
equivalent:
\begin{itemize}
\item[{\rm (i)}]
$ M-NR^+ N^T \geq 0 $, $R\geq0$, and  $ N(I-RR^+)=0$.
\item[{\rm (ii)}]
$\left[\ba{cc}
M & N \\
N^T & R
\ea\right] \geq 0$.
\item[\rm (iii)]
$\left[\ba{cc}
R & N^T \\
N & M
\ea\right] \geq 0$.
\end{itemize}
\end{lemma}

Its proof can be found in \cite{Alb:69}.

\ms

\bl\label{lem-equality} \sl
Let a constant matrix $M, N \in {\cal S}^n$ be given. Then for any admissible pair $(X(\cdot),u(\cdot))$
of the system (\ref{MF-state1}), we have
\bel{equa-M}\ba{ll}
\ns\ds \dbE\Big\{\int_0^t\Big[\big((X(s)-\dbE[X(s)]\big)^T(A^TM+MA+C^TMC)\big(X(s)-\dbE[X(s)]\big) \\
\ns\ds\q + 2\big(u(s)-\dbE[u(s)]\big)^T(B^TM+D^TMC)\big(X(s)-\dbE[X(s)]\big) \\
\ns\ds\q + \big(u(s)-\dbE[u(s)]\big)^TD^TMD\big(u(s)-\dbE[u(s)]\big) \\
\ns\ds\q + \big((C+\C)\dbE[X(s)]+(D+\D)\dbE[u(s)]\big)^TM\big((C+\C)\dbE[X(s)]+(D+\D)\dbE[u(s)]\big) \\
\ns\ds\q + \dbE[X(s)]^T\big((A+\A)^TN + N(A+\A)\dbE[X(s)]\big) + 2\dbE[X(s)]^T(A+\A)^TN(B+\B)\dbE[u(s)]\Big]ds\Big\} \\
\ns\ds = \dbE\Big[\big(X(t)-\dbE[X(t)]\big)^TM\big(X(t)-\dbE[x(t)]\big)\Big] + \dbE[X(t)]^TN\dbE[X(t)]-x^TNx, \;\;\forall t\geq 0. \ea\ee
\el

\it Proof. \rm Applying It\^o's formula to
$\big(X(t)-\dbE[X(t)]\big)^TM\big(X(t)-\dbE[X(t)]\big)$,
integrating from $0$ to $t$, and taking expectations, we easily get
the desired result.
\endpf

\ms

\bp\label{theo-stab} \sl
The following assertions are equivalent:
\begin{itemize}
\item[{\rm (i)}]
The controlled MF-FSDE system $[A,\bar A,C,\bar C;B,\bar B,D,\bar
D]$ is MF-$L^2$-stabilizable.
\item[{\rm (ii)}]
There exist matrices $K,\K$ and symmetric matrices $X,\X$ such that
\bel{cond1}\left\{\2n\ba{ll}
\ns\ds (A+BK)\mathbb{X} + \mathbb{X}(A+BK)^T + (C+D K)\mathbb{X}(C+D K)^T \\
\ns\ds\q + (C+\C+D\K+ \D\K)\bar{\mathbb{X}}(C+\C+D\K+ \D\K)^T < 0, \\
\ns\ds (A+\A+B\K+\B\K)\bar{\mathbb{X}} + \bar{\mathbb{X}}(A+\A+B\K+\B\K)^T < 0, \;\;\; \mathbb{X} > 0, \; \bar{\mathbb{X}} > 0.\ea\right.\ee
In this case the feedback $u(t)=K(X(t)-\dbE[X(t)])+\K\dbE[X(t)]$ is stabilizing.
\item[{\rm (iii)}]
There exist matrices $K,\K$ and symmetric matrices $\mathbb{X},\bar{\mathbb{X}}$ such that
\bel{cond2}\left\{\2n\ba{ll}
\ns\ds (A+BK)^T\mathbb{X} + \mathbb{X}(A+BK) + (C+D K)^T\mathbb{X}(C+D K) \\
\ns\ds\q + (C+\C+D\K+ \D\K)^T\bar{\mathbb{X}}(C+\C+D\K+ \D\K) < 0, \\
\ns\ds (A+\A+B\K+\B\K)\bar{\mathbb{X}} + \bar{\mathbb{X}}(A+\A+B\K+\B\K)^T < 0, \;\;\; \mathbb{X} > 0, \; \bar{\mathbb{X}} > 0. \ea\right.\ee
In this case the feedback $u(t)=K(X(t)-\dbE[X(t)])+\K\dbE[X(t)]$ is stabilizing.
\item[{\rm (iv)}]
There are matrices $K,\K$ such that for any matrices $Y,\Y$ there
exist unique solution $\mathbb{X},\bar{\mathbb{X}}$ to the following matrix equations
\bel{cond3}\left\{\2n\ba{ll}
\ns\ds (A+BK)\mathbb{X} + \mathbb{X}(A+BK)^T + (C+D K)\mathbb{X}(C+D K)^T \\
\ns\ds\q + (C+\C+D\K+ \D\K)\bar{\mathbb{X}}(C+\C+D\K+ \D\K)^T + Y = 0, \\
\ns\ds (A+\A+B\K+\B\K)\bar{\mathbb{X}} + \bar{\mathbb{X}}(A+\A+B\K+\B\K)^T + \Y = 0, \;\;\; \mathbb{X} > 0, \; \bar{\mathbb{X}} > 0. \ea\right.\ee
Moreover, if $Y, \Y > 0$ (resp. $Y, \Y \geq 0$) then $\mathbb{X},\bar{\mathbb{X}} > 0$ (resp. $\mathbb{X},\bar{\mathbb{X}} \geq 0$).
Furthermore, in this case the feedback $u(t)=K(X(t)-\dbE[X(t)])+\K\dbE[X(t)]$ is stabilizing.
\item[{\rm (v)}]
There are matrices $K,\K$ such that for any matrices $Y,\Y$ there
exist unique solution $\mathbb{X},\bar{\mathbb{X}}$ to the following matrix equations
\bel{cond4}\left\{\2n\ba{ll}
\ns\ds (A+BK)^T\mathbb{X} + \mathbb{X}(A+BK) + (C+D K)^T\mathbb{X}(C+D K) \\
\ns\ds\q + (C+\C+D\K+ \D\K)^T\bar{\mathbb{X}}(C+\C+D\K+ \D\K) + Y = 0, \\
\ns\ds (A+\A+B\K+\B\K)\bar{\mathbb{X}} + \bar{\mathbb{X}}(A+\A+B\K+\B\K)^T + \Y = 0, \;\;\; \mathbb{X} > 0, \; \bar{\mathbb{X}} > 0. \ea\right.\ee
Moreover, if $Y,\Y > 0$ (resp. $Y,\Y \geq 0$) then
$\mathbb{X},\bar{\mathbb{X}} > 0$ (resp. $\mathbb{X},\bar{\mathbb{X}} \geq 0$).
Furthermore, in this case the feedback $u(t)=K(X(t)-\dbE[X(t)])+\K\dbE[X(t)]$ is stabilizing.
\item[{\rm (vi)}]
There exist matrices $Y,\Y$ and symmetric matrices
$\mathbb{X},\bar{\mathbb{X}}$ such that
\bel{cond-LMI}\left\{\2n\ba{ll}\left[\ba{c|c}
\ns\ds \ba{ll} A\mathbb{X}+\mathbb{X}A^T+BY+Y^TB^T \\ +(C+\C+(D+\D)\Y\bar{\mathbb{X}}^{-1})\bar{\mathbb{X}}(C+\C+(D+\D)\Y\bar{\mathbb{X}}^{-1})^T \ea & C\mathbb{X}+DY \\ \hline
\mathbb{X}C^T+Y^TD^T & -\mathbb{X} \\ \ea\right] < 0, \\ [8mm]
\ns\ds (A+\A)\bar{\mathbb{X}}+(B+\B)\Y + \bar{\mathbb{X}}(A+\A)^T + \Y^T(B+\B)^T < 0, \qq \mathbb{X} > 0, \; \bar{\mathbb{X}} > 0. \ea\right.\ee
In this case the feedback
$u(t)=Y\mathbb{X}^{-1}(X(t)-\dbE[X(t)])+\Y\bar{\mathbb{X}}^{-1}\dbE[X(t)]$
is stabilizing.
\end{itemize}
\ep

\it Proof. \rm For any $n_u\times n$ matrices $K,\K$, define an
operator $\Phi, \widehat\Phi : {\cal S}^n\rightarrow {\cal S}^n$ by
$$\left\{\2n\ba{ll}
\ns\ds \Phi(\mathbb{X},\bar{\mathbb{X}}) = (A+BK)\mathbb{X} + \mathbb{X}(A+BK)^T + (C+D K)\mathbb{X}(C+D K)^T \\
\ns\ds\qq\qq\q + (C+\C+D\K+ \D\K)\bar{\mathbb{X}}(C+\C+D\K+ \D\K)^T, \\
\ns\ds \widehat\Phi(X,\bar{\mathbb{X}}) = (A+\A+B\K+\B\K)\bar{\mathbb{X}} + \bar{\mathbb{X}}(A+\A+B\K+\B\K)^T. \ea\right.$$
If  $X(\cdot)$ satisfies the equation (\ref{MF-state1}) with the
feedback control $u(t)=K(X(t)-\dbE[X(t)])+\K\dbE[X(t)]$, then by
It\^o's formula $\mathbb{X}(t) = \dbE\big[(X(t)-\dbE[X(t)])(X(t)-\dbE[X(t)])^T\big]$ and
$\bar{\mathbb{X}}(t) = \dbE[X(t)]\dbE[X(t)]^T$ satisfy the differential matrix systems
$$\frac{d}{dt}\mathbb{X}(t)=\Phi\big(\mathbb{X}(t),\bar{\mathbb{X}}(t)\big) \quad\mbox{ and }\quad \frac{d}{dt}\bar{\mathbb{X}}(t)=\widehat\Phi\big(\mathbb{X}(t),\bar{\mathbb{X}}(t)\big).$$
Applying the general result given in the appendix of \cite{ElA:96},
we have the equivalence between the mean-square stabilizability and
each of the assertions (ii)-(v). Furthermore, with $Y=K\mathbb{X}$
and $\Y=\K\bar{\mathbb{X}}$ the condition $(\ref{cond2})$ is
equivalent to
$$\left\{\2n\ba{ll}
\ns\ds A\mathbb{X}+\mathbb{X}A^T+BY+Y^TB^T + (C\mathbb{X}+DY)\mathbb{X}^{-1}(C\mathbb{X}+DY) \\
\ns\ds\q + (C+\C+(D+\D)\Y\bar{\mathbb{X}}^{-1})\bar{\mathbb{X}}(C+\C+(D+\D)\Y\bar{\mathbb{X}}^{-1})^T < 0,\\
\ns\ds (A+\A)\bar{\mathbb{X}}+(B+\B)\Y + \bar{\mathbb{X}}(A+\A)^T + \Y^T(B+\B)^T < 0, \qq \mathbb{X} > 0, \; \bar{\mathbb{X}} > 0. \ea\right.$$
Applying Schur's lemma (Lemma \ref{lem-schur})  we have the equivalence of the assertion (vi).
\endpf

\ms

Let $p^*$ denote the infimum value of the primal SDP (\ref{SDP-def})
and $d^*$ the supremum value of its dual (\ref{Dual-def}). Then we
have the following  results (\cite{VaB:96,AitZhou}).

\bp\label{theo-duality} \sl
$p^*=d^*$ if either of the following conditions holds:
\begin{itemize}
\item[{\rm (i)}]
The primal problem (\ref{SDP-def}) is strictly feasible, i.e., there exists an $x$ such that $F(x)> 0$.
\item[{\rm (ii)}]
The dual problem (\ref{Dual-def}) is strictly feasible, i.e., there exists a $Z\in {\cal S}^n$ with $Z> 0$ and
${\bf Tr}(ZF_i)=c_i,\;i=1,\cdots,m$.
\end{itemize}
If both conditions (i) and (ii) hold, then the optimal sets of both
the primal and the dual are nonempty. In this case, the following
complementary slackness condition
\bel{slackness-cond}
F(x)Z=0\ee
is necessary and sufficient for achieving the optimal values for both problems.
\ep

Now we turn to rewrite the AREs (\ref{A-Riccati}) as
\bel{ric-general}
\cR(P) = 0, \qq \bar\cR(P,\Pi) = 0, \ee
where
$$\left\{\2n\ba{ll}
\ns\ds \cR(P) \deq PA+A^TP+C^TPC-(PB+C^TPD)(R+D^TPD)^{-1}(B^TP+D^TPC)+Q, \\
\ns\ds \bar\cR(P,\Pi) \deq \Pi(A \1n+\1n \A) \1n+\1n (A \1n+\1n \A)^T\Pi \1n+\1n (C \1n+\1n \C)^TP(C \1n+\1n \C) \1n+\1n Q \1n+\1n \Q - \big[\Pi(B \1n+\1n \B) \1n+\1n (C \1n+\1n \C)^TP(D \1n+\1n \D)\big] \\
\ns\ds \qq\qq\q\cdot\big[R \1n+\1n \R\1n + \1n(D \1n+\1n\D)^TP(D \1n+\1n\D)\big]^{-1}\big[(B \1n+\1n \B)^T\Pi \1n+ \1n(D \1n+\1n \D)^TP(C\1n+\1n \C)\big]. \ea\right.$$
In this subsection, we pose an additional assumption that the interior of the set
$$\cP = \big\{(P,\Pi)\in {\cal S}^n\times{\cal S}^n \; \vert \; \cR(P)\geq 0, \bar\cR(P,\Pi) \geq 0\big\}$$
is nonempty, namely, there exists a $(P_0,\Pi_0)\in \cS^n\times{\cal S}^n$ such that $\cR(P_0)>0$, and $\bar\cR(P_0,\Pi_0) \geq 0$.

\ms

Consider the following SDP problem
\bel{general-LMI-Pb}\ba{rl}
\ns\ds \max & {\bf Tr}(P) + {\bf Tr}(\Pi), \\
\ns\ds \mbox{subject to} & \left\{\2n\ba{ll}\left[\ba{c|c} PA+A^TP+C^TPC+Q & PB+C^TPD \\ \hline B^TP+D^TPC & R+D^TPD \ea\right] \geq 0, \\
\ns\ds \left[\begin{array}{c|c}\ba{ll}\Pi(A + \A) + (A + \A)^T\Pi \\ + (C + \C)^TP(C + \C) + Q + \Q \ea & \Pi(B + \B)+(C + \C)^TP(D + \D) \\ \hline
(B + \B)^T\Pi+(D + \D)^TP(C + \C)  & R + \R + (D + \D)^TP(D + \D)\ea\right] \geq 0, \\
\ns\ds P-P_0 \geq 0, \\
\ns\ds \Pi-\Pi_0 \geq 0. \ea\right.\ea\ee

\ms

The constraints of SDP (\ref{general-LMI-Pb}) can be  equivalently
expressed as a single LMI
\bel{single-LMI}\ba{ll}
F(P,\Pi) \deq \left[\ba{c|c|c|c}
L(P) & 0 & 0 & 0 \\ \hline
0 & \L(P,\Pi) & 0 & 0 \\ \hline
0 & 0 & P-P_0 & 0 \\ \hline 0
& 0 & 0 & \Pi-\Pi_0 \ea\right] \geq 0, \ea\ee
where
$$\ba{ll}
\ns\ds L(P) \deq \left[\ba{c|c} PA+A^TP+C^TPC+Q & PB+C^TPD \\ \hline B^TP+D^TPC & R+D^TPD \ea\right], \\ [5mm]
\ns\ds \L(P,\Pi) \deq \left[\ba{c|c}\ba{ll}\Pi(A + \A) + (A + \A)^T\Pi \\ + (C + \C)^TP(C + \C) + Q + \Q \ea & \Pi(B + \B)+(C + \C)^TP(D + \D) \\ \hline
(B + \B)^T\Pi+(D + \D)^TP(C + \C)  & R+\R+(D+\D)^TP(D+\D) \ea\right]. \ea$$

\bp\label{dual} \sl
The dual problem of SDP (\ref{general-LMI-Pb}) can be formulated as follows
\bel{dual-Pb-general}\ba{rl}
\ns\ds \max & -{\bf Tr}(QS+WP_0+\Q\S+\W\Pi_0)-{\bf Tr}(RV+\R\V), \\
\ns\ds \mbox{\rm subject to} & \left\{\2n\ba{ll}\ba{ll}
\ns\ds AS+SA^T+BU+U^TB^T+CSC^T+DUC^T+CU^TD^T+DVD^T+(C+\C)\S(C+\C)^T \\
\ns\ds\q +(D+\D)\U(C+\C)^T+(C+\C)\U^T(D+\D)^T+(D+\D)\V(D+\D)^T+W+I = 0, \\
\ns\ds (A+\A)\S+\S(A+\A)^T+(B+\B)\U+\U^T(B+\B)^T+\W+I = 0, \ea \\
\ns\ds \left[\begin{array}{cc} S & U^T \\ U & V \ea\right] \geq 0, \;\;
\left[\ba{cc} \S & \U^T \\ \U & \V \ea\right] \geq 0, \;\;
W \geq 0, \;\; \W \geq 0, \ea\right.\ea\ee
where $S, \S, W, \W \in \cS^n$, $V, \V \in {\cal S}^m$ and $U, \U \in \dbR^{m\times n}$.
\ep

\it Proof. \rm The constraints of the general dual problem
(\ref{Dual-def}) can be formulated equivalently as the constraints
of (\ref{dual-Pb-general}). To this end, define the dual variable
$Z\in{\cal S}^{4n+2m}$ for (\ref{Dual-def}) as
$$Z=\left[\ba{c|c|c|c}
\begin{array}{cc}
S & U^T \\
U & T
\ea & Y_1^T & Y_2^T & Y_3^T \\ \hline
Y_1 & \ba{cc}
\S & \U^T \\
\U & \T
\ea & Y_4^T & Y_5^T \\ \hline
Y_2 & Y_4 & W & Y_6^T \\ \hline Y_3 & Y_5 & Y_6 & \W
\ea\right] \geq 0.$$
By the general duality relation ${\bf Tr}(ZF_i)=c_i,i=1,\cdots,m$
(see (\ref{Dual-def})) it follows that for any $(P,\Pi) \in \cS^n\times\cS^n$,
$${\bf Tr}([F(P,\Pi)-F(0,0)]Z)=-{\bf Tr}(P)-{\bf Tr}(\Pi),$$
which is equivalent to
$$\ba{ll}
\ns\ds {\bf Tr}\Big(\big[AS+SA^T+BU+U^TB^T+CSC^T+DUC^T+CU^TD^T+DVD^T+(C+\C)\S(C+\C)^T \\
\ns\ds\q +(D+\D)\U(C+\C)^T+(C+\C)\U^T(D+\D)^T+(D+\D)\V(D+\D)^T+W+I\big]P \\
\ns\ds\q +\big[(A+\A)\S+\S(A+\A)^T+(B+\B)\U+\U^T(B+\B)^T+\W+I\big]\Pi\Big) = 0. \ea$$
This leads to
$$\left\{\2n\ba{ll}
\ns\ds AS+SA^T+BU+U^TB^T+CSC^T+DUC^T+CU^TD^T+DVD^T+(C+\C)\S(C+\C)^T \\
\ns\ds\q +(D+\D)\U(C+\C)^T+(C+\C)\U^T(D+\D)^T+(D+\D)\V(D+\D)^T+W+I = 0, \\
\ns\ds (A+\A)\S+\S(A+\A)^T+(B+\B)\U+\U^T(B+\B)^T+\W+I = 0. \ea\right.$$
On the other hand, the objective of the dual problem
(\ref{Dual-def}) can be formulated as
$$-{\bf Tr}(F(0)Z)=-{\bf Tr}(QS+WP_0+\Q\S+\W\Pi_0)-{\bf Tr}(RV+\R\V).$$
In particular, since the matrix variables $Y_1, Y_2, Y_3, Y_4, Y_5$
and $Y_6$ do not play any role in the above formulation, they can be
dropped. Hence, the condition $Z\geq 0$ is equivalent to
$$\left[\ba{cc}
S & U^T\\
U & V
\ea\right] \geq 0, \;\;
\left[\ba{cc}
\S & \U^T\\
\U & \V
\ea\right] \geq 0, \;\;
W \geq 0, \;\; \W \geq 0.$$
This completes the proof.
\endpf

\ms

We now show that the MF-$L^2$-stability can be regarded as a dual
concept of SDP optimality.

\bp\label{theo-dual-strict} \sl
The dual problem (\ref{dual-Pb-general}) is strictly feasible if and only if the
controlled MF-FSDE system $[A,\bar A,C,\bar C;B,\bar B,D,\bar D]$ is
MF-$L^2$-stabilizable.
\ep

\it Proof. \rm First, assume that the controlled MF-FSDE system
$[A,\bar A,C,\bar C;B,\bar B,D,\bar D]$ is MF-$L^2$-stabilizable by
some feedback $u(t)=K\big(X(t)-\dbE[X(t)]\big)+\K\dbE[X(t)]$. Let
$\widetilde{W} > 0$ and $\widehat{W} > 0$ be fixed. Then it follows
from the assertion (v) of Proposition \ref{theo-stab} that there
exists a unique $(S,\S)$ satisfying
$$\left\{\2n\ba{ll}
\ns\ds (A+BK)S + S(A+BK)^T + (C+D K)S(C+D K)^T \\
\ns\ds\q + (C+\C+D\K+ \D\K)\S(C+\C+D\K+ \D\K)^T + \widetilde{W} + I = 0, \\
\ns\ds (A+\A+B\K+\B\K)\S + \S(A+\A+B\K+\B\K)^T + \widehat{W} + I = 0, \q S > 0, \; \S > 0. \ea\right.$$
Set $U=KS$ and $\U=\K\S$. The above relation can then be rewritten
as
$$\left\{\2n\ba{ll}
\ns\ds AS + SA + BU + U^TB^T + CSC^T + DUC^T + CU^TD^T + DUS^{-1}U^TD^T + (C+\C)\S(C+\C)^T \\
\ns\ds\q + (D+\D)\U(C+\C)^T+(C+\C)\U^T(D+\D)^T+(D+\D)\U\S^{-1}\U^T(D+\D)^T + \widetilde{W} + I = 0, \\
\ns\ds (A+\A)\S + \S(A+\A)^T +(B+\B)\U+\U^T(B+\B)^T + \widehat{W} + I = 0. \ea\right.$$
Let  $\epsilon > 0$ and $\bar\epsilon > 0$, define
$V=\epsilon I+US^{-1}U^T$, $\V=\bar\epsilon I+\U\S^{-1}\U^T$,
$W= -\epsilon DD^T-\bar\epsilon(D+\D)(D+\D)^T+\widetilde{W}$ and
$\widehat{W} = \W$. Then $V$, $\V$, $W$ and $\W$ satisfy
$$\left\{\2n\ba{ll}
\ns\ds AS + SA + BU + U^TB^T + CSC^T + DUC^T + CU^TD^T + DVD^T  + (C+\C)\S(C+\C)^T \\
\ns\ds\q + (D+\D)\U(C+\C)^T+(C+\C)\U^T(D+\D)^T+(D+\D)\V(D+\D)^T + W + I = 0, \\
\ns\ds (A+\A)\S + \S(A+\A)^T +(B+\B)\U+\U^T(B+\B)^T + \W + I = 0. \ea\right.$$
Moreover, by Schur's lemma (Lemma \ref{lem-schur}) for $\epsilon >0$
and $\bar\epsilon > 0$ sufficiently small we must have
$$\left[\ba{cc}
S & U^T \\
U & V
\ea\right] \geq 0, \;\;
\left[\ba{cc}
\S & \U^T \\
\U & \V
\ea\right] \geq 0, \;\;
W \geq 0, \;\; \W \geq 0.$$
Therefore, the dual problem (\ref{dual-Pb-general}) is strictly feasible.

Conversely, assume that the dual problem is strictly feasible. Then
there exist $S>0, \S>0, U, \U, V$ and $\V$ such that
$$\left\{\ba{ll}
\ns\ds AS + SA + BU + U^TB^T + CSC^T + DUC^T + CU^TD^T + DVD^T  + (C+\C)\S(C+\C)^T \\
\ns\ds\q + (D+\D)\U(C+\C)^T+(C+\C)\U^T(D+\D)^T+(D+\D)\V(D+\D)^T < 0, \\
\ns\ds (A+\A)\S + \S(A+\A)^T+(B+\B)\U+\U^T(B+\B)^T < 0. \ea\right.$$
It follows that
$$\left\{\ba{ll}
\ns\ds AS + SA + BU + U^TB^T + CSC^T + DUC^T + CU^TD^T + DUS^{-1}U^TD^T  + (C+\C)\S(C+\C)^T \\
\ns\ds\q + (D+\D)\U(C+\C)^T+(C+\C)\U^T(D+\D)^T+(D+\D)\U\S^{-1}\U^T(D+\D)^T < 0, \\
\ns\ds (A+\A)\S+ \S(A+\A)^T+(B+\B)\U+\U^T(B+\B)^T < 0. \ea\right.$$
Define $K=US^{-1}$ and $\K=\U\S^{-1}$. The above inequality is equivalent to
$$\left\{\begin{array}{ll}
\ns\ds (A+BK)S + S(A+BK)^T + (C+D K)S(C+D K)^T \\
\ns\ds\q + (C+\C+D\K+ \D\K)\S(C+\C+D\K+ \D\K)^T < 0, \\
\ns\ds (A+\A+B\K+\B\K)\S + \S(A+\A+B\K+\B\K)^T < 0, \q S > 0, \; \S > 0.
\end{array}\right.$$
We conclude that the assertion (iii) of Proposition \ref{theo-stab}
is satisfied. Hence, the controlled MF-FSDE system $[A,\bar A,C,\bar
C;B,\bar B,D,\bar D]$ is MF-$L^2$-stabilizable.
\endpf

\ms

The following result presents the existence of the solution of the
AREs (\ref{ric-general}) via the SDP (\ref{general-LMI-Pb}).

\bp\label{theo-sol1} \sl
The optimal set of SDP
(\ref{general-LMI-Pb}) is nonempty and any optimal solution
$(P_*,\Pi_*)$ must satisfy the ARE (\ref{ric-general}).
\ep

\it Proof. \rm Proposition \ref{theo-dual-strict}, along with
Proposition \ref{theo-duality}, yields the non-emptiness of the
optimal set. Next, appealing to the complementary slackness
condition (\ref{slackness-cond}) in Proposition \ref{theo-duality},
we conclude that any optimal solution $(P_*,\Pi_*)$ must satisfy
$$\left[\ba{c|c|c|c}
L(P) & 0 & 0 & 0 \\ \hline
0 & \L(P,\Pi) & 0 & 0 \\ \hline
0 & 0 & P-P_0 & 0 \\ \hline
0 & 0 & 0 & \Pi-\Pi_0\ea\right]
\left[\ba{c|c|c|c}\ba{cc}
S & U^T \\ U & V \ea & Y_1^T & Y_2^T & Y_3^T \\ \hline
Y_1 & \ba{cc} \S & \U^T \\ \U & \V \ea & Y_4^T & Y_5^T \\ \hline
Y_2 & Y_4 & W & Y_6^T \\ \hline
Y_3 & Y_5 & Y_6 & \W \ea\right]= 0,$$
where $S,\S,U,\U,V,\V,W$ and $\W$ are the corresponding optimal dual
variables. From the above we can deduce the following conditions
\bel{opt-cond1}
(A^TP_*+P_*A+C^TP_*C+Q)S+(P_*B+C^TP_*D)U=0, \ee
\bel{opt-cond2}
(A^TP_*+P_*A+C^TP_*C+Q)U^T+(P_*B+C^TP_*D)V=0, \ee
\bel{opt-cond3}
(B^TP_*+D^TP_*C)S+(R+D^TP_*D)U=0,
\ee
\bel{opt-cond4}
(B^TP_*+D^TP_*C)U^T+(R+D^TP_*D)V=0,
\ee
\bel{opt-cond5}
\Big(\Pi_*(A + \A) + (A + \A)^T\Pi_* + (C + \C)^TP_*(C + \C) + Q + \Q\Big)\S+\Big(\Pi(B + \B)+(C + \C)^TP_*(D + \D)\Big)\U=0, \ee
\bel{opt-cond6}
\Big(\Pi_*(A + \A) + (A + \A)^T\Pi_* + (C + \C)^TP_*(C + \C) + Q + \Q\Big)\U^T+\Big(\Pi(B + \B)+(C + \C)^TP_*(D + \D)\Big)\V=0, \ee
\bel{opt-cond7}
\Big((B + \B)^T\Pi_*+(D + \D)^TP_*(C + \C)\Big)\S + \Big(R+\R+(D+\D)^TP_*(D+\D)\Big)\U=0, \ee
\bel{opt-cond8}
\Big((B + \B)^T\Pi_*+(D + \D)^TP_*(C + \C)\Big)\U^T + \Big(R+\R+(D+\D)^TP_*(D+\D)\Big)\V=0, \ee
\bel{opt-cond9}
(P_*-P_0)W=0, \ee
\bel{opt-cond10}
(\Pi_*-\Pi_0)\W=0. \ee
%
Hence (\ref{opt-cond3}) implies that
$U=-(R+D^TP_*D)^{-1}(B^TP_*+D^TP_*C)S$. Putting this into equation
(\ref{opt-cond1}) leads to $\cR(P_*)S=0$. A same
manipulation of equations (\ref{opt-cond2}) and (\ref{opt-cond4})
yields $\cR(P_*)U^T=0$.
Similarly, (\ref{opt-cond7}) implies that
$\U=-\big(R+\R+(D+\D)^TP_*(D+\D)\big)^{-1}\big[(B + \B)^T\Pi_*+(D + \D)^TP_*(C + \C)\big]\S$. Substituting this into equation
(\ref{opt-cond5}) leads to $\bar\cR(P_*,\Pi_*)\S=0$. And a similar
same manipulation of equations (\ref{opt-cond6}) and
(\ref{opt-cond8}) yields $\bar\cR(P_*,\Pi_*)\U^T=0$. Recall that the
dual variables $S,\S,U,\U,V,\V,W,\W$ satisfy the following constraint
\bel{equa-proof-theo-stab-sol-1}\ba{ll}
\ns\ds AS + SA + BU + U^TB^T + CSC^T + DUC^T + CU^TD^T + DVD^T  + (C+\C)\S(C+\C)^T \\
\ns\ds\q + (D+\D)\U(C+\C)^T+(C+\C)\U^T(D+\D)^T+(D+\D)\V(D+\D)^T + W + I = 0. \ea\ee
Multiplying both sides  of the above by $\cR(P_*,\Pi_*)$ we have
$$\ba{ll}
\ns\ds \cR(P_*)[CSC^T + DUC^T + CU^TD^T + DVD^T  + (C+\C)\S(C+\C)^T + (D+\D)\U(C+\C)^T \\
\ns\ds \q+(C+\C)\U^T(D+\D)^T+(D+\D)\V(D+\D)^T + W + I]\cR(P_*)=0. \ea$$
It follows from $W \geq 0$ that
\bel{equa-proof}\ba{ll}
\ns\ds \cR(P_*)[CSC^T + DUC^T + CU^TD^T + DUS^+\U^TD^T  + (C+\C)\S(C+\C)^T + (D+\D)\U(C+\C)^T \\
\ns\ds\q +(C+\C)\U^T(D+\D)^T+(D+\D)\U\S^+\U^T(D+\D)^T]\cR(P_*) \leq 0. \ea\ee
Since
$$\left[\ba{cc}
S & U^T \\
U & T
\ea\right] \geq 0,
\quad \left[\ba{cc}
\S & \U^T \\
\U & \T
\ea\right] \geq 0,$$
it follows from extended Schur's lemma (Lemma \ref{schur-ext}) that
$V\geq US^+U^T$, $U=USS^+$, $\V\geq \U\S^+\U^T$ and $\U=\U\S\S^+$.
By virtue of Lemma \ref{pseudo-lem} we deduce the following
\bel{inequa-proof-1}\ba{ll}
\ns\ds CSC^T+DUC^T+CU^TD^T+DUS^+ U^TD^T \\
\ns\ds = CSS^+ SC^T+DUS^+ SC^T+CSS^+ U^TD^T+ DUS^+ U^TD^T \\
\ns\ds = (CS+DU)S^+ (SC^T+U^TD^T) \geq 0, \ea\ee
and
\bel{inequa-proof-2}\ba{ll}
\ns\ds (C+\C)\S(C+\C)^T + (D+\D)\U(C+\C)^T + (C+\C)\U^T(D+\D)^T+(D+\D)\U\S^+\U^T(D+\D)^T \\
\ns\ds = (C+\C)\S\S^+\S(C+\C)^T+(D+\D)\U\S^+\S(C+\C)^T+(C+\C)\S\S^+\U^T(D+\D)^T \\
\ns\ds\q + (D+\D)\U\S^+\U^T(D+\D)^T \\
\ns\ds = \Big((C+\C)\S+(D+\D)\U\Big)\S^+\Big(\S(C+\C)^T+\U^T(D+\D)^T\Big) \geq 0. \ea\ee
Then it follows from (\ref{equa-proof}) that $\cR(P_*)\cR(P_*)\leq 0$, resulting in $\cR(P_*)=0$.

Recall that the dual variables $S,\S,T,\T,U,\U,W,\W$ satisfy the
following constraint
\begin{equation}\label{equa-proof-theo-stab-sol-2}
(A+\A)\S + \S(A+\A)^T+(B+\B)\U+\U^T(B+\B)^T + \W + I = 0.
\end{equation}
Multiplying   both sides  of the above by $\bar\cR(P_*,\Pi_*)$, we
have
$$\begin{array}{l}
\bar\cR(P_*,\Pi_*)[\W + I]\bar\cR(P_*,\Pi_*)=0.
\end{array}$$
Since $\W \geq 0$, we have $\bar\cR(P_*,\Pi_*)=0$.
\endpf

\ms

The following result indicates that any optimal  solution of the
primal SDP gives rise to an MF-$L^2$ stabilizing control of the
MF-LQ problem. The readers can refer to \cite{BDDM}.

\ms

\bp\label{theo-stab-sol} \sl
Let $(P_*,\Pi_*)$ be an optimal
solution to the primal  SDP (\ref{general-LMI-Pb}). Then the
feedback control $u(t)=\Gamma_*\big(X(t)-\dbE[X(t)]\big)+\bar\Gamma_*X(t)$ is
stabilizing for the system (\ref{MF-state1}), where
$$\left\{\2n\ba{ll}
\ns\ds \Gamma_* = -(R+D^TP_*D)^{-1}(B^TP_*+D^TP_*C), \\
\ns\ds \bar\Gamma_*=-\big(R+\R+(D+\D)^TP_*(D+\D)\big)^{-1}\big((B + \B)^T\Pi_*+(D + \D)^TP_*(C + \C)\big). \ea\right.$$
\ep

\it Proof. \rm Let $S,\S,U,\U,V,\V,W$ and $\W$ be the corresponding
optimal dual variables satisfying
(\ref{opt-cond1})-(\ref{opt-cond10}). First, we are to show that
$S>0$ and $\S>0$. Suppose that $Sx=0$ and $\S x=0, x\in \dbR^n$. As
$U$ and $\U$ satisfy
$$U=-(R+D^TP_*D)^{-1}(B^TP_*+D^TP_*C)S$$
and
$$\U=-\big(R+\R+(D+\D)^TP_*(D+\D)\big)^{-1}\big[(B + \B)^T\Pi_*+(D + \D)^TP_*(C + \C)\big]$$
(see (\ref{opt-cond3}) and (\ref{opt-cond7})), we also have  $Ux=0$
and $\U x=0$. The dual constraint (\ref{equa-proof-theo-stab-sol-1})
then implies
$$\ba{ll}
\ns\ds x^T\big[CSC^T + DUC^T + CU^TD^T + DUS^+\U^TD^T  + (C+\C)\S(C+\C)^T + (D+\D)\U(C+\C)^T \\
\ns\ds\q +(C+\C)\U^T(D+\D)^T+(D+\D)\U\S^+\U^T(D+\D)^T+W+I\big]x \leq 0. \ea$$
The same manipulation as in the proof of Proposition \ref{theo-sol1}
gives $x=0$. As  $S\geq 0$ and $\S \geq 0$, we conclude that $S>0$
and $\S>0$. Now, the equalities (\ref{equa-proof-theo-stab-sol-1})
and (\ref{equa-proof-theo-stab-sol-2}) give
$$\left\{\2n\ba{ll}
\ns\ds AS + SA + BU + U^TB^T + CSC^T + DUC^T + CU^TD^T + DUS^{-1}U^TD^T  + (C+\C)\S(C+\C)^T \\
\ns\ds\q + (D+\D)\U(C+\C)^T+(C+\C)\U^T(D+\D)^T+(D+\D)\U\S^{-1}\U^T(D+\D)^T < 0, \\
\ns\ds (A+\A)\S + \S(A+\A)^T+(B+\B)\U+\U^T(B+\B)^T < 0, \q S > 0, \;\; \S > 0, \ea\right.$$
which is equivalent to the mean-square stabilizability condition
({iii}) of Proposition  \ref{theo-stab} with $K=\Gamma_*$ and $\K =
\bar\Gamma_*$.
\endpf

\bp\label{coro-sol1} \sl
There exists a unique optimal solution to
the SDP (\ref{general-LMI-Pb}), which is also the maximal solution
to the AREs (\ref{ric-general}).
\ep

\it Proof. \rm Let $(P_*,\Pi_*)$ be an optimal solution to the SDP
(\ref{general-LMI-Pb}). Proposition \ref{theo-sol1} shows that
$(P_*,\Pi_*)$ solves the AREs (\ref{ric-general}). To show that it
is indeed a maximal solution, define
$$\left\{\2n\ba{ll}
\ns\ds \Gamma_* = -(R+D^TP_*D)^{-1}(B^TP_*+D^TP_*C), \\
\ns\ds \bar\Gamma_* = -\big(R+\R+(D+\D)^TP_*(D+\D)\big)^{-1}\big[(B + \B)^T\Pi_*+(D + \D)^TP_*(C + \C)\big]. \ea\right.$$
A simple calculation yields
$$\left\{\2n\ba{ll}
\ns\ds (A+B\Gamma_*)^TP_* + P_*(A+B\Gamma_*) + (C+D\Gamma)^TP_*(C+D\Gamma_*) = -Q-\Gamma_*^TR\Gamma_*, \\
\ns\ds (A+\A+B\bar\Gamma_* + \B\bar\Gamma_*)\Pi_* + \Pi_*(A+\A+B\bar\Gamma_*+\B\bar\Gamma_*)^T \\
\ns\ds\q + (C+\C+D\bar\Gamma_* + \D\bar\Gamma_*)^TP_*(C+\C+D\bar\Gamma_*+\D\bar\Gamma_*) = -Q-\Q-\bar\Gamma_*^T(R+\R)\bar\Gamma_*. \ea\right.$$
On the other hand, it follows from Proposition \ref{theo-stab-sol}
that
$u_*(t)=\Gamma_*\big(X_*(t)-\dbE[X_*(t)]\big)+\bar\Gamma_*\dbE[X_*(t)]$
is a stabilizing control. A proof similar to that of Theorem
\ref{theo-vcost=Pmax} yields that $(P_*,\Pi_*)$ is the upper bound
of the set ${\cal P}$, namely, $(P_*,\Pi_*)$ is the maximal
solution. Finally, the uniqueness of the solution to the SDP
(\ref{general-LMI-Pb}) follows from the maximality.
\endpf

\end{document}